\documentclass[12pt]{amsart}
\usepackage{amsmath}
\usepackage{amsxtra}
\usepackage{amscd}
\usepackage{amsthm}
\usepackage{amsfonts}
\usepackage{amssymb}
\usepackage{eucal}

\input prepictex
\input pictex
\input postpictex

\newtheorem{thm}{Theorem}[section]
\newtheorem{lem}{Lemma}[section]
\newtheorem{cor}{Corollary}[section]
\newtheorem{prop}{Proposition}[section]

\theoremstyle{definition}

\newtheorem{example}{Example}[section]

\theoremstyle{remark}
\newtheorem{rem}{Remark}[section]

\numberwithin{equation}{section}

\begin{document}

\newcommand{\thmref}[1]{Theorem~\ref{#1}}
\newcommand{\secref}[1]{Section~\ref{#1}}
\newcommand{\lemref}[1]{Lemma~\ref{#1}}
\newcommand{\propref}[1]{Proposition~\ref{#1}}
\newcommand{\corref}[1]{Corollary~\ref{#1}}
\newcommand{\remref}[1]{Remark~\ref{#1}}
\newcommand{\eqnref}[1]{(\ref{#1})}
\newcommand{\exref}[1]{Example~\ref{#1}}

\newcommand{\nc}{\newcommand}
\nc{\on}{\operatorname} \nc{\Z}{{\mathbb Z}} \nc{\C}{{\mathbb C}}
\nc{\oo}{{\mf O}}

\nc{\bib}{\bibitem} \nc{\pa}{\partial} \nc{\F}{{\mf F}}
\nc{\rarr}{\rightarrow} \nc{\larr}{\longrightarrow}
\nc{\al}{\alpha} \nc{\ri}{\rangle} \nc{\lef}{\langle}
\nc{\W}{{\mc W}} \nc{\gam}{\ol{\gamma}} \nc{\Q}{\ol{Q}}
\nc{\q}{\widetilde{Q}} \nc{\la}{\lambda} \nc{\ep}{\epsilon}
\nc{\g}{\mf g} \nc{\h}{\mf h}
\nc{\n}{\mf n} \nc{\A}{{\mf a}}
\nc{\G}{{\mf g}} \nc{\Li}{{\mc L}} \nc{\La}{\Lambda}
\nc{\is}{{\mathbf i}} \nc{\V}{\mf V} \nc{\bi}{\bibitem}
\nc{\NS}{\mf N} \nc{\dt}{\mathord{\hbox{${\frac{d}{d t}}$}}}
\nc{\E}{\mc E} \nc{\ba}{\tilde{\pa}}
\def\smapdown#1{\big\downarrow\rlap{$\vcenter{\hbox{$\scriptstyle#1$}}$}}
\nc{\mc}{\mathcal} \nc{\mf}{\mathfrak} \nc{\ol}{\fracline}
\nc{\el}{\ell} \nc{\etabf}{{\bf \eta}} \nc{\x}{{\bf x}}
\nc{\xibf}{{\bf \xi}} \nc{\y}{{\bf y}} \nc{\half}{{\frac{1}{2}}}
\advance\headheight by 2pt \nc{\SN}{\mf S}
\nc{\ov}[1]{\overline{#1}}

\title{Finite conformal modules over {$N=2,3,4$} superconformal algebras}

\author[Shun-Jen Cheng]{Shun-Jen Cheng}
\address{Department of Mathematics, National Cheng Kung University,
Tainan, Taiwan}\email{chengsj@mail.ncku.edu.tw}

\author{Ngau Lam}
\address{Department of Mathematics, National Cheng Kung University, Tainan,
Taiwan} \email{nlam@mail.ncku.edu.tw}

\thanks{both authors acknowledge partial support by NSC-grant 90-2115-M-006-002 of the R.O.C}

\begin{abstract}
In this paper we continue the study of representation theory of
formal distribution Lie superalgebras initiated in \cite{CK1}. We
study finite Verma-type conformal modules over the $N=2$, $N=3$
and the two $N=4$ superconformal algebras and also find
explicitly all singular vectors in these modules.  From our
analysis of these modules we obtain a complete list of finite
irreducible conformal modules over the $N=2$, $N=3$ and the two
$N=4$ superconformal algebras. \vspace{.3cm}

\noindent{\bf Mathematics Subject Classifications (1991)}: 17B65
\vspace{.3cm}

\noindent{\bf PACS}: 11.10.Kk; 11.25.Hf

\vspace{.3cm}

\noindent{\bf Key words}: Superconformal algebras,
representations, conformal modules

\end{abstract}

\maketitle

\section{Introduction}

Superconformal algebras  have been playing an important role in
the study of string theory and conformal field theory, which have
been the subject of intensive study since the seminal paper
\cite{BPZ}. Superconformal algebras may be viewed as natural
super-extensions of the Virasoro algebra and their roots in
physics literature can be traced at least back to as early as the
70's \cite{A}. A mathematically rigorous definition of a
superconformal algebra is as follows. It is a simple Lie
superalgebra $\G$ over the complex numbers $\C$ spanned by the
modes of a finite family $\F$ of mutually local fields satisfying
the following two axioms \cite{K1}:
\begin{itemize}
\item[1.] $\F$ contains the Virasoro field,
\item[2.] the coefficients of the operator product expansions of members
from $\F$ are linear combinations of members from $\F$ and their
derivatives.
\end{itemize}
A Lie superalgebra $\G$ satisfying the second axiom only is
referred to as a \emph{formal distribution Lie superalgebra} in
\cite{K1}.

In order to facilitate the study of formal distribution Lie
superalgebras the notion of a \emph{conformal superalgebra} was
introduced in \cite{K1} (see \secref{prelim}).  It proves to be an
effective tool for this purpose.

A natural class of representations of formal distribution Lie
superalgebras to study is the class of \emph{conformal modules}
\cite{CK1}.  A conformal module is a pair consisting of a
$\G$-module $V$ and a family $\E$ of fields whose modes span $V$
such that members from $\F$ and $\E$ are mutually local. Just as
the study of formal distribution Lie superalgebras reduces to the
study of conformal superalgebras, the study of conformal modules
is essentially reduced to the study of modules over the
corresponding conformal superalgebras.

The study of modules over the conformal superalgebra can further
be reduced to the study of modules over the \emph{extended
annihilation subalgebra}, which is a semidirect sum of the
subalgebra of positive modes of the corresponding formal
distribution Lie superalgebra and a one-dimensional derivation.
It is in this language that the problem of classifying finite
irreducible conformal modules over the Virasoro, $N=1$
(Neveu-Schwarz) and the current superalgebra was solved in
\cite{CK1}.

The problem of classifying conformal modules over other
superconformal algebras, which is the main theme of the present
paper, turns out to be more subtle. The main purpose here is to
give a classification of finite irreducible conformal modules
over the $N=2$, $N=3$ and the two $N=4$ superconformal algebras.

We first construct finite Verma-type conformal modules for a
general superconformal algebra and prove that every finite
irreducible conformal module is a homomorphic image of such a
module. As a consequence we obtain a bijection between finite
irreducible conformal modules of a superconformal algebra and
finite-dimensional irreducible modules of a certain
finite-dimensional reductive Lie (super)algebra
(\corref{bijection}).

We then study these Verma-type modules in detail for the four
members of the family of superconformal algebras mentioned above.
It turns out that, unlike for the Virasoro and the $N=1$
(Neveu-Schwarz) superconformal algebras, the Verma-type modules
for these superconformal algebras are in general reducible, and
thus we need to analyze their submodules. This is accomplished by
finding explicit formulas for all singular vectors inside such a
module and then show that the submodule generated by these
singular vectors is maximal (in all but two cases). We also find
an explicit basis for this maximal submodule, which then enables
us to give a quite explicit description of all finite irreducible
conformal modules over these superconformal algebras.

This paper is organized as follows. In \secref{prelim} basic
facts of formal distribution Lie superalgebras, conformal
superalgebras and extended annihilation subalgebras are recalled.
\secref{mainlemma} is devoted to the study of a class of modules
over a certain class of Lie superalgebras that include the
annihilation subalgebra of every superconformal algebra. This
class of modules gives rise to finite Verma-type conformal modules
of superconformal algebras.  The results of \secref{mainlemma}
are then used in \secref{N=2}, \secref{N=3}, \secref{N=4} and
\secref{BN=4}, where finite irreducible conformal modules over
the $N=2$, $N=3$, the ``small'' $N=4$ and the ``big'' $N=4$
superconformal algebra, respectively, are classified.

In this paper all vector spaces, (super)algebras and tensor
products are over taken over the complex numbers $\C$.

\section{Preliminaries}\label{prelim}

In this section we review some of the basic facts on formal
distribution Lie (super)algebras and conformal modules that will
be used later on.  The material here is taken from \cite{CK1},
\cite{K1} and \cite{K3}, and the reader is referred to these
articles for more details.

\subsection{Formal Distribution Lie
Superalgebras}\label{conformalalgebra}

Recall that a {\it formal distribution} or a \emph{field} with
coefficients in a Lie superalgebra $\G=\G_{\bar{0}}+\G_{\bar{1}}$
is a formal series of the form:
$$a(z)=\sum_{n\in\Z}a_{[n]}z^{-n-1},$$ where $a_{[n]}\in \G$ and
$z$ is an indeterminate.

Two formal distributions $a(z)$ and $b(z)$ with coefficients in
$\G$ are said to be mutually \emph{local} if there exists
$N\in\Z_+$ such that
\begin{equation}\label{orieqn}
(z-w)^N[a(z),b(w)]=0.
\end{equation}

Let $\delta(z-w)=z^{-1}\sum_{n\in\Z}(\frac{z}{w})^n$ be the formal
delta function. Then \eqnref{orieqn} may be written as
\begin{equation}\label{OPEeqn}
[a(z),b(w)]=\sum_{j=0}^{N-1}(a_{(j)}b)(w){{\pa_{w}^{(j)}}\delta(z-w)},
\end{equation}
(here $\pa^{(j)}_w$ stands for
${\frac{1}{j!}}{\frac{\partial^j}{\partial w^j}}$) for some
uniquely determined formal distributions $(a_{(j)}b)(w)$, and thus
defines a $\C$-bilinear product $\cdot_{(j)}\cdot$ for each
$j\in\Z_+$ on the space of all formal distributions with
coefficients in $\G$.  Also $\pa_z a(z)=\sum_n(\pa
a)_{[n]}z^{-n-1}$, where $(\pa a)_{[n]}=-na_{[n-1]}$, and hence
the space of all formal distributions is also a (left)
$\C[\pa_z]$-module.

A Lie superalgebra $\G$ is called a {\it formal distribution Lie
superalgebra}, if there exists a family $\F$ of mutually local formal
distributions whose coefficients span $\G$. We will write $(\G,\F)$
for such a Lie superalgebra.

Given a formal distribution Lie superalgebra $(\G,\F)$, we may
include $\F$ in the minimal family $\overline{\F}$ of mutally
local distributions which is closed under $\pa_{z}$ and all
products $\cdot_{(j)}\cdot$.  Then $\overline{\F}$ is a {\it
conformal superalgebra}, i.e.~it is a left $\Z_2$-graded
$\C[\pa]$-module $R$ with a $\C$-bilinear product $a_{(n)}b$ for
each $n\in\Z_+$ such that the following axioms hold ($a,b,c\in R;
m,n\in\Z_+$ and $\pa^{(j)}={\frac{1}{j!}}\pa^j$) (cf.~\cite{B},
\cite{FLM}):
\begin{itemize}
\item[(C0)] $a_{(n)}b=0$, for $n>>0$,
\item[(C1)] $(\pa a)_{(n)}b=-na_{(n-1)}b$,
\item[(C2)] $a_{(n)}b=(-1)^{p(a)p(b)}\sum_{j=0}^{\infty}(-1)^{j+n+1}{\pa^{(j)}}(b_{(n+j)}a)$,
\item[(C3)] $a_{(m)}(b_{(n)}c)=\sum_{j=0}^{\infty}{\binom{m}
{j}}(a_{(j)}b)_{(m+n-j)}c+(-1)^{p(a)p(b)}b_{(n)}(a_{(m)}c)$.
\end{itemize}
It is convenient to write the products of $a,b\in R$ in the
generating series form $$a_\la b=\sum_{n=0}^\infty
a_{(n)}b\frac{\la^n}{n!},$$ where $\la$ is a formal
indeterminate.  Such an expression lies in $R[\la]$.

Conversely, if a conformal superalgebra $R=\bigoplus_{i\in
I}\C[\pa]a^i$ is free $\C[\pa]$-module, we may associate to $R$ a
formal distribution Lie superalgebra $(\G(R),\F(R))$ with Lie
superalgebra $\G(R)$ spanned by $\C$-basis $a^i_{[m]}$ ($i\in I$,
$m\in\Z$) and fields
$\F(R)=\{a^i(z)=\sum_{n\in\Z}a^i_{[n]}z^{-n-1}\}_{i\in I}$ with
bracket (cf.~\eqref{OPEeqn}):
\begin{equation*}
[a^i(z),a^j(w)]=\sum_{k\in\Z_+}(a^i_{(k)}a^j)(w){{\pa_w^{(k)}}}\delta(z-w),
\end{equation*}
so that $\overline{\F(R)}=R$, giving rise to  commutation
relations ($m,n\in\Z$; $i,j\in I$)
\begin{equation}\label{generic1}
[a^i_{[m]},a^j_{[n]}]=\sum_{k\in\Z_+}{\binom{m}
{k}}(a^i_{(k)}a^j)_{[m+n-k]}.
\end{equation}
It follows that the Lie superalgebra $\G$ of a formal
distribution Lie superalgebra $(\G,\F)$ is isomorphic to
$\G(\overline{\F})$ divided by an \emph{irregular} ideal, that is
an ideal which does not contain every $a_{[n]}$ for some non-zero
element $a\in\overline{\F}$.

\begin{example} The (centerless) {\it Virasoro
algebra} $\V$ has a basis $L_n$ ($n\in\Z$) with commutation
relations $$[L_m,L_n]=(m-n)L_{m+n}.$$ It is spanned by the
coefficients of the field $L(z)=\sum_{n\in\Z}L_n z^{-n-2}$
satisfying
\begin{equation}\label{VirOPE}
[L(z),L(w)]=\pa_{w}L(w)\delta(z-w)+2L(w)\pa_{w}\delta(z-w).
\end{equation}
The conformal algebra associated to the Virasoro algebra, is the
{\it Virasoro conformal algebra} $R(\V)=\C[\pa]\otimes L$ with
products $L_\la L=(\pa +2\la)L$.
\end{example}

\begin{example}
Let $\G$ be a finite-dimensional Lie (super)algebra.  Let
$\tilde{\G}=\G\otimes\C[t,t^{-1}]$ denote the corresponding
\emph{current algebra} with bracket
\begin{equation*}
[a\otimes f(t),b\otimes g(t)]=[a,b]\otimes f(t)g(t), \quad
a,b\in\G; f(t),g(t)\in\C[t,t^{-1}].
\end{equation*}
For each $a\in\G$ define a field $a(z)=\sum_{n\in\Z}(a\otimes
t^{n})z^{-n-1}$.  Then $\tilde{\G}$ is spanned by the coefficients
of $a(z)$ satisfying
\begin{equation}\label{currentOPE}
[a(z),b(w)]=[a,b](w)\delta(z-w).
\end{equation}
The conformal (super)algebra associated to the current algebra is
the {\it current conformal algebra} $R(\tilde{\G})=\C[\pa]\otimes
\G$ with products $a_\la b=[a,b]$, $a,b\in\G$.
\end{example}

\begin{example}
The semidirect sum $\V\ltimes\tilde{\G}$ is another example of a
formal distribution Lie (super)algebra.  The collection of fields
is $\{L(z),a(z)|a\in\G\}$ and we have in addition to
\eqnref{VirOPE} and \eqnref{currentOPE}
\begin{equation}\label{semiOPE}
[L(z),a(w)]=\pa_wa(w)\delta(z-w)+a(w)\pa_w\delta(z-w).
\end{equation}
The conformal algebra associated to the semidirect sum of the
Virasoro algebra and the current algebra is
$R(\V\ltimes\tilde{\G})=R(\V)\ltimes R(\tilde{\G})$. For $a\in\G$
we have $L_\la a=(\pa +\la)a$.
\end{example}

\subsection{Conformal Modules}\label{conformalmodule}

Let $(\G,\F)$ be a formal distribution Lie superalgebra.  Let $V$
be a $\G$-module such that $V$ is spanned over $\C$ by the
coefficients of a family $\E$ of fields.  If all $a(z)\in\F$ are
local with respect to all $v(z)\in \E$, then the pair $(V,\E)$ is
called a {\it conformal module} over $(\G,\F)$.

Now the family $\E$ of a conformal module $(V,\E)$ over $(\G,\F)$
similarly can be included in a larger family $\overline{\E}$,
which is still local with respect to the fields from
$\overline{\F}$, and invariant under $\pa$ and $a_{(j)}$, for all
$a\in\overline{\F}$ and $j\in\Z_+$. It can be shown that for
$a,b\in\overline{\F}$ and $v\in\overline{\E}$ ($m,n\in\Z_+$) one
has $$[a_{(m)},b_{(n)}]v=\sum_{j=0}^m
{\binom{m}{j}}(a_{(j)}b)_{(m+n-j)}v,\quad(\pa a)_{(n)}
v=[\pa,a_{(n)}]v=-na_{(n-1)}v.$$ Thus it follows that any
conformal module $(V,\E)$ over a formal distribution Lie
superalgebra $(\G,\F)$ gives rise to a module $M=\overline{\E}$
over the \emph{conformal superalgebra} $R=\overline{\F}$, defined
as follows. It is a (left) $\Z_2$-graded $\C[\pa]$-module equipped
with a family of $\C$-linear maps $a\rightarrow a_{(n)}^M$ of $R$
to ${\rm End}_{\C}M$, for each $n\in\Z_+$, such that the following
properties hold for $a,b\in R$ and $m,n\in\Z_+$:
\begin{itemize}
\item[(M0)] $\quad a^M_{(n)}v=0$, for $v\in M$ and $n>>0$,
\item[(M1)] $\quad [a^M_{(m)},b^M_{(n)}]=\sum_{j=0}^m{\binom{m}
{j}}(a_{(j)}b)^M_{(m+n-j)}$,
\item[(M2)] $\quad (\pa a)^M_{(n)}=[\pa,a^M_{(n)}]=-na_{(n-1)}^M$.
\end{itemize}
Again it is convenient to write the action of an element $a\in R$
on an element $v\in M$ in the form of a generating series in
$V[\la]$
\begin{equation*}
a_\la v:=\sum_{n=0}^{\infty} a_{(n)}v \frac{\la^n}{n!}.
\end{equation*}

Conversely, suppose that a conformal superalgebra
$R=\bigoplus_{i\in I}\C[\pa]a^i$ is a free $\C[\pa]$-module and
consider the associated formal distribution Lie superalgebra
$(\G(R),\F(R))$. Let $M$ be a module over the conformal
superalgebra $R$ and suppose that $M$ is a free $\C[\pa]$-module
with $\C[\pa]$-basis $\{v^{{\alpha}}\}_{{\alpha}\in J}$. This
gives rise to a conformal module $V(M)$ over $\G(R)$ with fields
$\E=\{v^\alpha(z)=\sum_{n\in\Z}v^\alpha_{[n]}z^{-n-1}|\alpha\in
J\}$ and $\C$-basis $v_{[n]}^{{\alpha}}$, defined by:
$$a^{i}(z)v^{{\alpha}}(w)=\sum_{j\in\Z_+}(a^{i}_{(j)}v^{{\alpha}})(w){\pa_w^{(j)}}\delta(z-w).$$

A conformal module $(V,\E)$ (respectively module $M$) over a
formal distribution Lie superalgebra $(\G,\F)$ (respectively over
a conformal superalgebra $R$) is called {\it finite}, if
$\overline{\E}$ (respectively $M$) is a finitely generated
$\C[\pa]$-module. A conformal module $(V,\E)$ over $(\G,\F)$ is
called {\it irreducible}, if there is no non-trivial invariant
subspace which contains all $v_{[n]}$, $n\in\Z$, for some non-zero
$v\in \overline{\E}$.  An invariant subspace that does not
contain all $v_{[n]}$, for some non-zero $v\in\E$, is called an
\emph{irregular submodule} and conformal modules that differ by
an irregular submodule are called referred to as
\emph{equivalent} in \cite{K3}. Clearly a conformal module is
irreducible if and only if the associated module $\overline{\E}$
over the conformal superalgebra $\overline{\F}$ is irreducible.

\begin{rem} \label{invariant} It follows from (M2) that an eigenvector
$v\in M$ of the linear operator $\pa$ is an $R$-invariant,
i.e.~$a_{(n)}v=0$, for all $n\ge 0$. Thus a finite irreducible
module over a conformal superalgebra $R$ is either free over
$\C[\pa]$ or else it is one-dimensional over $\C$.
\end{rem}

Suppose that $(\G,\F)$ is a formal distribution Lie superalgebra
such that $\G(\overline{\F})\cong\G$. Our discussion implies that
any irreducible conformal module $(V,\E)$ over $(\G,\F)$ is a
quotient of an irreducible conformal module of the form $V(M)$
divided by an irregular submodule, where $M$ is an irreducible
module over the conformal superalgebra $\overline{\F}$. Hence in
particular if $V(M)$ is irreducible as a $\G$-module for every
irreducible $M$, then every finite irreducible conformal modules
over $(\G,\F)$ isomorphic to $V(M)$, for some finite irreducible
$\overline{\F}$-module $M$.

\begin{example}\label{virconf} The Virasoro algebra $\V$ may be identified with
the Lie algebra of regular vector fields on $\C^{\times}$, where
$L_n=-t^{n+1}\dt$, $n\in\Z$.  For $\alpha,\Delta\in\C$ let
$$F_{\V}(\alpha,\Delta)=\C[t,t^{-1}]e^{-\alpha t}dt^{1-\Delta}.$$
The Lie algebra $\V$ acts on the space $F_{\V}(\alpha,\Delta)$ in
a natural way:
\begin{equation*}
(f(t)\frac{\pa}{\pa t}) g(t)dt^{1-\Delta}=(f(t)g'(t)+(1-\Delta)
g(t)f'(t))dt^{1-\Delta},
\end{equation*}
where $f(t)\in\C[t,t^{-1}]$ and $g(t)\in \C[t,t^{-1}]e^{-\alpha
t}$.  Letting $v_{[n]}=t^ne^{-\alpha t}dt^{1-\Delta}$ and
$v(z)=\sum_{n\in\Z}v_{[n]}z^{-n-1}$ this action is equivalent to
$$L(z)v(w)=(\pa_{w}+\alpha)v(w)\delta(z-w)+\Delta
v(w)\pa_{w}\delta(z-w).$$ Hence we have constructed a
two-parameter family of conformal modules over $\V$. This gives a
family of $R(\V)$-modules $\C[\pa]\otimes\C v_{\Delta}$ with
products $L_\la v_{\Delta}=(\alpha+\pa+\Delta\la)v_{\Delta}$.
This module is irreducible if and only if $\Delta\not=0$, in
which case it will be denoted by $L_{\V}(\alpha,\Delta)$.  We set
$L_{\V}(\alpha,0)$ to be the one-dimensional (over $\C$)
$R(\V)$-module on which $\pa$ acts as the scalar $\alpha$.
\end{example}

\begin{example}\label{currentconf}
Let $\G$ be a finite-dimensional simple Lie algebra and $U^\La$
the finite-dimensional irreducible module of highest weight
$\La$. Then $F_{\tilde{\G}}(\La)=U^\La\otimes\C[t,t^{-1}]$ is
naturally a module over $\tilde{\G}$ with action given by
\begin{equation}\label{sec2generic1}
(a\otimes f(t))(u\otimes g(t))=au\otimes f(t)g(t),\quad
a\in\G,u\in U^\La;f(t),g(t)\in\C[t,t^{-1}].
\end{equation}
For each vector $u\in U^\La$ define $u(z)=\sum_{n\in\Z}(u\otimes
t^n) z^{-n-1}$ so that \eqnref{sec2generic1} is equivalent to
\begin{equation*}
a(z)u(w)=au(w)\delta(z-w),
\end{equation*}
and hence $F_{\tilde{\G}}(\La)$ is conformal. This gives a family
of $R(\tilde{\G})$-modules, which is irreducible if and only if
$\La\not=0$, in which case it will be denoted by
$L_{\tilde{\G}}(\Lambda)$. By $L_{\tilde{\G}}(0)$ we will mean
the trivial $R(\tilde{\G})$-module. Similarly one defines the
one-dimensional module $L_{\tilde{\G}}(\alpha,0)$.
\end{example}

\begin{example}\label{semiconf}
$\tilde{\G}$ acts on
$F_{\V\ltimes\tilde{\G}}(\alpha,\Delta,\La)=U^\La\otimes
F_{\V}(\alpha,\Delta)$ similarly as in \exref{currentconf}.
However, on $F_{\V\ltimes\tilde{\G}}(\alpha,\Delta,\La)$ we have
also an action of $\V$, thus making it into a module over
$\V\ltimes \tilde{\G}$.  This module defines an
$R(\V\ltimes\tilde{\G})$-module which is irreducible if and only
if $(\Delta,\La)\not=(0,0)$, and in which case it will be denoted
by $L_{\V\ltimes\tilde{\G}}(\alpha,\Delta,\Lambda)$. By
$L_{\V\ltimes\tilde{\G}}(\alpha,0,0)$ we will mean the
one-dimensional module on which $\pa$ acts a the scalar $\alpha$.
\end{example}

The following theorem was proved in \cite{CK1}.

\begin{thm}\label{classvir} Let $\G$ stand for a
finite-dimensional simple Lie algebra. Any finite irreducible
module over the conformal algebras $R(\V)$, $R(\tilde{\G})$ and
$R(\V\ltimes\tilde{\G})$ are as follows:
\begin{itemize}
\item[i.] $L_{\V}(\alpha,\Delta)$,
\item[ii.] $L_{\tilde{\G}}(\La)$ and $L_{\tilde{\G}}(\alpha,0)$,
\item[iii.] $L_{\V\ltimes\tilde{\G}}(\alpha,\Delta,\La)$.
\end{itemize}
\end{thm}

\begin{rem}
We note that a similar statement as \thmref{classvir} part (iii)
holds even if $\G$ is replaced by the $1$-dimensional Lie algebra
$\C a$.  In this case $U^\La=\C u$ with $a u=\La u$, $\La\in\C$.
Also part (ii) remains true for all but three series of
finite-dimensional simple Lie superalgebras.
\end{rem}

\subsection{Extended Annihilation Subalgebras}\label{annihilation}
Given a formal distribution Lie superalgebra $(\G,\F)$ we let
$\G_+$ denote the $\C$-span of all $a_{[n]}$, where $n\ge 0$ and
$a\in\F$.  Due to \eqnref{generic1} $\G_+$ is closed under the
bracket and hence form a subalgebra of $\G$, which we will call
the {\it annihilation algebra} of $(\G,\F)$.  Let $\pa$ be the
derivation of $\G_+$ defined by $[\pa,a_{[n]}]=-na_{[n-1]}$, and
consider the semi-direct sum of $\G^+=\C\pa\ltimes \G_+$.  Then
$\G^+$ is called the \emph{extended annihilated algebra} of
$(\G,\F)$. The following proposition, which follows by comparing
(M1) with \eqnref{generic1}, is important for the theory of
conformal modules.

\begin{prop}\cite{CK1} Let $R$ be a conformal superalgebra and
$(\G(R),R(\F))$ be its associated formal distribution Lie
superalgebra with extended annihilation algebra $\G(R)^+$.  Then
a module over the conformal superalgebra $R$ is precisely a
$\G(R)^+$-module $M$ satisfying $a_{[n]}v=0$, for each $v\in M$,
$a\in R$ and $ n>>0$.
\end{prop}

\begin{rem}
Let $R$ be a conformal superalgebra with $\C[\pa]$-basis
$\{a^i|i\in I\}$ and $M$ a free $\C[\pa]$-module with basis
$\{v^j|j\in J\}$. Given $a^i_{(n)}v^j\in M$ for all $i\in I$,
$j\in J$, $n\in\Z_+$, which is $0$ for $n>>0$, condition (M2)
uniquely extends the action of $a^i_{(n)}$ to all of $M$.  If in
addition (M1) holds, then $M$ is an $R$-module.  Hence the action
of an $R$-module $M$ is completely determined by the action of a
$\C[\pa]$-basis of $R$ on a $\C[\pa]$-basis of $M$.
\end{rem}

\begin{example}\label{annihilationmodules}
In the case of the Virasoro algebra $\V$ the annihilation algebra
$\V_+$ is spanned by elements $L_n$, $n\ge -1$.  In the case of
the current algebra $\tilde{\G}_+$ is spanned by $a\otimes t^n$,
where $a\in\G$ and $n\ge 0$, while in the case of
$\V\ltimes\tilde{\G}$ it is $\V_+\ltimes\tilde{\G}_+$.
\end{example}

The problem of classifying conformal modules over $(\G,\F)$ is
thus reduced to the problem of classifying a class of modules over
$\G(\overline{\F})^+$.  It is clear that in all our examples one
has $\G(\overline{\F})=\G$, and thus we are to study modules over
$\G^+$.  Now if in addition there exists an element $L_{-1}$ in
$\G_+$ such that $L_{-1}-\pa$ is central in $\G^+$, then every
irreducible representation of $\G^+$ is an irreducible
representation of $\G_+$, on which $(L_{-1}-\pa)$ acts as a
scalar $\alpha\in\C$. In the case of the $\V$ and
$\V\ltimes\tilde{\G}$ and the $N=2,3,4$ superconformal
superalgebras, which we will define later, such an $L_{-1}$ always
exists so that we only need to consider representations of $\G_+$.
The irreducible representations of $\V_+$, and
$\V_+\ltimes\tilde{\G}_+$ that give rise to those in
\thmref{classvir} are denoted by $L_{\V_+}(\Delta)$ and
$L_{\V_+\ltimes\tilde{\G}_+}(\Delta,\La)$, respectively. The
corresponding actions are clear and can be found in \cite{CK1}.

\section{Finite Verma-type Conformal Modules}\label{mainlemma}

Let $\Li$ be a Lie superalgebra over $\C$ with a distinguished
element $\pa$ and a descending sequence of subspaces
$\Li=\Li_{-1}\supset\Li_{0}\supset\Li_{1}\supset\Li_{2}\supset\cdots\supset\Li_{n}
\supset\cdots$, such that $[\pa,\Li_{k}]=\Li_{k-1}$, for all
$k>0$. Let $W$ be an $\Li$-module, which is finitely generated
over $\C[\pa]$, such that for all $w\in W$ there exists a
non-negative integer $k$ (depending on $w$) with $\Li_{k}w=0$.
For $m\ge -2$ set $W_m=\{w\in W|\Li_{m+1}w=0\}$ and let $M$ be the
minimal non-negative integer such that $W_M\not=0$.

\begin{lem}\cite{CK1}\label{keylemma}
Suppose that $M\ge 0$.  Then $\C[\pa]W_M=\C[\pa]\otimes W_M$ and
hence $\C[\pa]W_M\cap W_M=W_M$.  In particular $W_M$ is a
finite-dimensional vector space.
\end{lem}

Let $\G$ be a Lie superalgebra satisfying the following three
conditions.
\begin{itemize}
\item[(L1)] $\G$ is $\Z$-graded of finite depth
$d\in\mathbb{N}$, i.e.~ $\G=\bigoplus_{j\ge -d}\G_{j}$ with
$[\G_i,\G_j]\subset\G_{i+j}$.
\item[(L2)] There exists a semisimple element
$z\in\G_0$ such that it centralizer in $\G$ is contained in
$\G_0$.
\item[(L3)] There exits an element $\pa\in\G_{-d}$ such that
$[\pa,\G_i]=\G_{i-d}$, for $i\ge 0$.
\end{itemize}

\begin{rem}
If $\G$ contains the grading operator with respect to its
gradation, then condition (L2) is automatic.
\end{rem}

Examples of Lie superalgebras satisfying (L1)--(L3) are provided
by annihilation subalgebras of superconformal algebras, which we
will describe in more detail.

Let $t$ be an even indeterminate and $\xi_1,\ldots,\xi_N$ be $N$
odd indeterminate.  Denote by $\Lambda(N)$ the Grassmann
superalgebra in the indeterminates $\xi_1,\ldots,\xi_N$ and set
$\Lambda(1,N):=\C[t,t^{-1}]\otimes\Lambda(N)$.  Let $W(1,N)$ be
the derivation superalgebra of $\Lambda(1,N)$, then $W(1,N)$ is a
formal distribution Lie superalgebra \cite{K2}.  Letting
$\frac{\pa}{\pa t}$ and $\frac{\pa}{\pa\xi_i}$, for
$i=1,\ldots,N$, be the usual differential operators, every element
in $D\in W(1,N)$ can be written as \cite{K4}
\begin{equation*}
D=a_0\frac{\pa}{\pa t}+\sum_{i=1}^N a_i\frac{\pa}{\pa\xi_i},\quad
a_0,a_i,\ldots,a_N\in\Lambda(1,N).
\end{equation*}
The \emph{standard gradation} of $W(1,N)$ is obtained by setting
the degree of $t$ and $\xi_i$ to be $1$. Its annihilation
subalgebra is $W(1,N)_+=\bigoplus_{j\ge -1}(W(1,N))_j$.
$W(1,N)_+$ in this gradation contains its grading operator given
by $z=t\frac{\pa}{\pa t}+\sum_{i=1}^N\xi_i\frac{\pa}{\pa\xi_i}$ so
that (L2) is satisfied. Also choosing $\pa$ to be $\frac{\pa}{\pa
t}$ it follows that (L3) is also satisfied so that $W(1,N)$ is a
Lie superalgebra of the type above.  Note that $W(1,N)_0\cong
gl(1,N)$.

The subalgebra of divergence zero vector fields in $W(1,N)$
contains an ideal of codimension $1$.  This ideal is its derived
algebra and is the superconformal algebra $S(1,N)$ \cite{K2}. The
standard gradation of $W(1,N)_+$ induces a gradation on the
annihilation subalgebra $S(1,N)_+$ of $S(1,N)$.  Choosing
$z=t\frac{\pa}{\pa
t}+\frac{1}{N}\sum_{i=1}^N\xi_i\frac{\pa}{\pa\xi_i}$ along with
$\pa=\frac{\pa}{\pa t}$ it follows that $S(1,N)_+$ in this
gradation also satisfies (L1)--L(3). Observe that $S(1,N)_0\cong
sl(1,N)$ and also that the ``small'' $N=4$ superconformal algebra
(to be defined in \secref{N=4}) is isomorphic to $S(1,2)$
\cite{KL}.

The contact superalgebra $K(1,N)$ is the subalgebra of $W(1,N)$
defined by
\begin{equation*}
K(1,N):=\{D\in W(1,N)|D\omega=f_D\omega,\ {\rm for\ some\
}f_D\in\Lambda(1,N)\},
\end{equation*}
where $\omega:=dt-\sum_{i=1}^N\xi_i d\xi_i$ is the standard
contact form. Here the action of $D$ on $\omega$ is the usual
action of vector fields on differential forms.

The map from $\Lambda(1,N)$ to $K(1,N)$ given by to
\begin{equation*}
f\rightarrow 2f\frac{\pa}{\pa
t}+(-1)^{p(f)}\sum_{i=1}^N(\xi_i\frac{\pa f}{\pa t}+\frac{\pa
f}{\pa\xi_i})(\xi_i\frac{\pa }{\pa t}+\frac{\pa }{\pa\xi_i})
\end{equation*}
is a bijection and hence it allows us to identify $K(1,N)$ with
the polynomial superalgebra $\Lambda(1,N)$.  The Lie bracket in
$\Lambda(1,N)$, also called the contact bracket, then reads for
homogeneous elements $f,g\in\Lambda(1,N)$:
\begin{equation*}
[f,g]=(2-E)f\frac{\pa g}{\pa t}-\frac{\pa f}{\pa
t}(2-E)g+(-1)^{p(f)}\sum_{i=1}^N\frac{\pa f}{\pa\xi_i}\frac{\pa
g}{\pa\xi_i},
\end{equation*}
where $E=\sum_{i=1}^N\xi_i\frac{\pa}{\pa\xi_i}$ is the Euler
operator.

When $N$ is even it is sometimes more convenient to make the
change of basis
$\xi^+_j=\frac{1}{\sqrt{2}}(\xi_j+i\xi_{j+\frac{N}{2}})$ and
$\xi^-_j=\frac{1}{\sqrt{2}}(\xi_j-i\xi_{j+\frac{N}{2}})$, for
$j=1,\ldots,\frac{N}{2}$ and $i=\sqrt{-1},$ so that the contact
bracket takes the split form:
\begin{equation*}
[f,g]=(2-E)f\frac{\pa g}{\pa t}-\frac{\pa f}{\pa t}(2-E)g
+(-1)^{p(f)}\sum_{i=1}^{\frac{N}{2}}(\frac{\pa
f}{\pa\xi^+_i}\frac{\pa g}{\pa\xi^-_i}+\frac{\pa
f}{\pa\xi^-_i}\frac{\pa g}{\pa\xi^+_i}),
\end{equation*}
where $E$ again is the Euler operator
$\sum_{i=1}^{\frac{N}{2}}(\xi^+_i\frac{\pa}{\pa\xi^+_i}
+\xi^-_i\frac{\pa}{\pa\xi^-_i})$.

The contact superalgebra $K(1,N)$ is a formal distribution Lie
superalgebra with fields defined as follows: Let
$I=\{i_1,\ldots,i_k\}$ be an ordered subset of $\{1,\ldots,N\}$,
and denote by $\xi_I$ the monomial $\xi_{i_1}\cdots\xi_{i_k}$.
Each such monomial gives rise to a field
$\xi_I(z)=\sum_{j\in\Z}\xi_It^j z^{-j-1}$. Evidently the span of
the coefficients of all such $\xi_I(z)$ is $K(1,N)$.  Furthermore
it is easy to check that these fields are mutually local and form
a formal distribution Lie superalgebra. This Lie superalgebra
becomes $\Z$-graded by putting the degree of $\xi_I t^n$ to
$2n+k-2$.  Obviously $t$ is the grading operator of this
gradation.  This gradation of $K(1,N)$ is usually referred to as
its \emph{standard gradation}.

The annihilation subalgebra $K(1,N)_+$ of $K(1,N)$ is spanned by
the basis elements $\xi_I t^n$, where $n\ge 0$ and $I$ runs over
all subsets of $\{i_1,\ldots,i_k\}$ ordered in (strictly)
increasing order. The $\Z$-gradation from $K(1,N)$ induces a
gradation on $K(1,N)_+$ making it a $\Z$-graded Lie superalgebra
of depth $2$ so that
$K(1,N)_+=\bigoplus_{j=-2}^{\infty}(K(1,N)_+)_j$ satisfies (L1)
and (L2).  In this gradation it is easy to check that
$[1,(K(1,N)_+)_j]=(K(1,N)_+)_{j-2}$ for all $j\ge 0$, so that
$K(1,N)_+$ also satisfies condition (L3).  It is easy to see that
the annihilation subalgebra of the small $N=4$ superconformal
algebra, which we define in \secref{N=4}, also satisfies
conditions (L1)--(L3). Note that $K(1,N)_0\cong cso_N$, the
direct sum of the Lie algebra $so_N$ and the one-dimensional Lie
algebra.

Finally it follows from the description of the exceptional
superconformal algebra $CK_6$ as a subalgebra of $K(1,6)$ in
\cite{CK2} that its annihilation subalgebra
$(CK_6)_+=\bigoplus_{j\ge -2} (CK_6)_j$ is a Lie superalgebra
satisfying (L1)--(L3) with $(CK_6)_0\cong cso_6$.

The modules over the annihilation subalgebras that are equivalent
to modules over the corresponding conformal superalgebras are then
$\G$-modules $V$ satisfying the following conditions.
\begin{itemize}
\item[(V1)] For all $v\in V$ there
exists an integer $k_0\ge -d$ (depending on $v$) such that
$\G_{k}v=0$, for all $k\ge k_0$.
\item[(V2)] $V$ is finitely generated over $\C[\pa]$.
\end{itemize}

We shall call $\G$-modules satisfying these two properties
\emph{finite}. Let $V$ be a finite irreducible $\G$-module. For
$n\ge -d-1$ set $V_n=\{v\in V|\G_jv=0,\forall j>n\}$. Let $N$ be
the minimal integer such that $V_N\not=0$. Such an $N$ exists by
(V1).

\begin{lem} \label{finite-dim} If $N\ge 0$, then
$V_N$ is a finite-dimensional vector space over $\C$.
\end{lem}

\begin{proof}
We let $\Li=\G$ and put $\Li_j=\bigoplus_{i\ge jd}\G_i$ so that we
have a filtration of subspaces
$$\Li\supset\Li_{0}\supset\Li_{1}\supset\Li_{2}\supset\cdots\supset\Li_{n}
\supset\cdots,$$ with $[\pa,\Li_{i}]=\Li_{i-1}$, for all $i\ge 0$
by (L3).  Let $W_m:=\{v\in V|\Li_{m+1} v=0\}$ and let $M$ be the
minimal integer such that $W_M\not=0$.  Since $N\ge 0$ implies
that $M\ge 0$, this setting puts us in the situation of
\lemref{keylemma}, from which we conclude that $W_M$ is a
finite-dimensional vector space over $\C$.  Of course $V_N\subset
W_M$ and hence it follows that $V_N$ is finite-dimensional as
well.
\end{proof}

We obtain the following description of finite irreducible
$\G$-modules.

\begin{thm}\label{key}
Let $\G=\bigoplus_{j\ge -d}\G_j$ be a Lie superalgebra satisfying
conditions (L1)--(L3) and $V$ a finite irreducible $\G$-module.
There exists a finite-dimensional irreducible $\G_0$-module
$U_0$, extended trivially to an $\Li_0(={\bigoplus_{j\ge
0}\G_j})$-module, and a $\G$-epimorphism $\varphi:{\rm
Ind}^\G_{\Li_0} U_0\rightarrow V$.
\end{thm}

\begin{proof}
We will continue to use the notation defined earlier. First we
show that $N\le 0$. Suppose that $N>0$. It is easy to see that
$V_N$ is invariant under $\Li_0$. Now there exits a basis
$\{x_1,\ldots,x_m\}$ of $\G_N$ together with non-zero complex
number $\la_1,\ldots,\la_m$ such that $[z,x_i]=\la_i x$, where $z$
is the element of (L2). Since $V_N$ is a finite-dimensional vector
space it follows in particular that $x_i$ acts nilpotently on
$V_N$ for all $1\le i\le m$. But
$[\G_N,\G_N]\subset\bigoplus_{j\ge N+1}\G_{j}$ and so the action
of the $x_i$'s on $V_N$ commutes. Therefore there exits a non-zero
$v\in V_N$ such that $\G_N v=0$.  But in this case
$V_{N-1}\not=0$, which contradicts the minimality of $N$. Thus
$N\le 0$.

In the case when $N=0$, there exists an epimorphism of
$\G$-modules ${\rm Ind}_{\Li_0}^{\G}V_0\rightarrow V$, with $V_0$
finite-dimensional due to \lemref{finite-dim}. By irreducibility
of $V$ it follows that $V_0=U_0$ is an irreducible $\G_0$-module.
Now if $N<0$, then there exists a non-zero vector $v$ invariant
under the action of $\G_j$, for $j\ge 0$. Again we have an
epimorphism of $\G$-modules ${\rm Ind}_{\Li_0}^{\G}\C
v\rightarrow V$.
\end{proof}

As a corollary of \thmref{key} we obtain the following.

\begin{cor}\label{bijection} There exists a bijection between
finite irreducible conformal modules of the superconformal
algebra $\G$ and finite-dimensional irreducible representations
of the Lie (super)algebra $\G_0$, where
\begin{itemize}
\item[i.] $\G=K(1,N)$ and $\G_0=cso_N$,
\item[ii.] $\G=W(1,N)$ and $\G_0=gl(1,N)$,
\item[iii.] $\G=S(1,N)$ and $\G_0=sl(1,N)$,
\item[iv.] $\G=CK_6$ and $\G_0=cso_6$.
\end{itemize}
\end{cor}

\begin{proof}
By \thmref{key} every finite irreducible $\G$-module is a
homomorphic image of ${\rm Ind}_{\Li_0}^{\G}U_0$. Now the usual
argument for highest weight representations implies that given a
finite-dimensional irreducible $\G_0$-module $U_0$ the
$\G$-module ${\rm Ind}_{\Li_0}^{\G}U_0$ contains a unique maximal
submodule, from which the bijection then follows.
\end{proof}

\begin{rem}
It is usual to put a half-integer gradation on $K(1,N)$ when
thinking of it as a superconformal algebra. The grading operator
of $K(1,N)$ with respect to this gradation is then $\frac{t}{2}$
rather than $t$. In this gradation one has
$K(1,N)_+=\bigoplus_{j\ge -1}\G_j$, where $j\in\half\Z$.
\thmref{key} of course remains valid after making some obvious
changes regarding gradation.  For a Lie superalgebra
$\G=\bigoplus_{j\ge -1}\G_j$ with $j\in\half\Z$, we will make it a
convention to write $\G_-$ for the subalgebra
$\bigoplus_{j<0}\G_j$.
\end{rem}

\section{Finite irreducible Modules over the $N=2$ conformal
superalgebra}\label{N=2}

The $N=2$ superconformal algebra is the formal distribution Lie
superalgebra $K(1,2)$.  Letting $\xi^+,\xi^-$ denote the two odd
indeterminates (so that we are using the split contact form) this
algebra is generated by the following four fields:
$L(z)=\sum_{n\in\Z}-\frac{t^{n+1}}{2} z^{-n-2}$,
$G^\pm(z)=\sum_{r\in\frac{1}{2}+\Z}\xi^\pm t^{r+\frac{1}{2}}
z^{-r-\frac{3}{2}}$ and $J(z)=\sum_{n\in\Z}\xi^-\xi^+t^n
z^{-n-1}$.  Its corresponding conformal superalgebra is then
generated freely over $\C[\pa]$ by $\{L,J,G^{\pm}\}$ with
products:
\begin{align*}
&L_\la L=(\pa+2\la)L,\quad L_\la J=(\pa+\la)J,\quad L_\la
G^{\pm}=(\pa+\frac{3}{2}\la)G^{\pm},\allowdisplaybreaks\\
&J_\la G^{\pm}=\pm G^{\pm},\quad G^+_\la G^-=(\pa+2\la)J+2L.
\end{align*}
Letting $L_{n}=-\frac{t^{n+1}}{2}$, $G^\pm_r=\xi^\pm
t^{r+\frac{1}{2}}$ and $J_n=\xi^-\xi^+ t^n$ with $n\in\Z$, $r\in
\frac{1}{2}+\Z$, the non-zero brackets in $K(1,2)$ are
($m,n\in\Z$ and $r,s\in\frac{1}{2}+\Z$):
\begin{align*}
&[L_m,L_n]=(m-n)L_{m+n},\quad
[L_m,G^{\pm}_r]=(\frac{m}{2}-r)G^{\pm}_{m+r},\quad [L_m,J_n]=-nJ_{n+m},\\
&[J_m,G^{\pm}_r]={\pm}G^{\pm}_{m+r},\qquad
[G^+_r,G^-_s]=2L_{r+s}+(r-s)J_{r+s}.
\end{align*}

The annihilation subalgebra $\G=K(1,2)_+$ is then spanned by
$L_{m}$, $J_n$ and $G^{\pm}_r$, where $m\ge -1$, $n\ge 0$ and
$r\ge -\frac{1}{2}$.  Note that letting $\G_j$ be the span of
$X_j$, where $X=L,J,G^{\pm}$, equips $\G=\bigoplus_{j\ge -1}\G_j$,
$j\in\frac{1}{2}\Z$, with a (consistent)
$\frac{1}{2}\Z$-gradation.  We denote $L_{-1}$ by $\pa$ from now
on.

Let $\C v_{\Delta,\La}$, $\Delta,\La\in\C$, be the one-dimensional
module over the abelian Lie algebra $\G_0=\C L_0+\C J_0$,
determined by
\begin{equation*}
L_0 v_{\Delta,\La}=\Delta v_{\Delta,\La},\quad J_0
v_{\Delta,\La}=\La v_{\Delta,\La}.
\end{equation*}
We may extend $\C v_{\Delta,\La}$ to a module over
$\Li_0=\bigoplus_{j\ge 0}\G_j$ by setting $\G_j v_{\Delta,\La}=0$,
for $j>0$.  Let $M_{\NS^2_+}(\Delta,\La):={\rm Ind}_{\Li_0}^{\G}\C
v_{\Delta,\La}$.  We denote by $N$ the unique maximal submodule of
$M_{\NS^2_+}(\Delta,\Lambda)$. The quotient
$M_{\NS^2_+}(\Delta,\Lambda)/N$ is the irreducible highest weight
module $L_{\NS^2_+}(\Delta,\Lambda)$ of highest weight
$(\Delta,\Lambda)$.  By \thmref{key} $L_{\NS^2_+}(\Delta,\Lambda)$
for $\Delta,\Lambda\in\C$ form a complete list of finite
irreducible $K(1,2)_+$-modules. Our next objective is to give a
more explicit description of $N$ and hence of
$L_{\NS^2_+}(\Delta,\Lambda)$.

It is clear that $\pa^k v_{\Delta,\La}$, $\pa^k
G^+_{-\frac{1}{2}}v_{\Delta,\La}$, $\pa^k
G^-_{-\frac{1}{2}}v_{\Delta,\La}$ and $\pa^k
G^+_{-\frac{1}{2}}G^-_{-\frac{1}{2}}v_{\Delta,\La}$, $k\ge 0$, is
a basis consisting of $(L_0,J_0)$-weight vectors for
$M_{\NS^2_+}(\Delta,\La)$ of $(L_0,J_0)$-weights $(\Delta+k,\La)$,
$(\Delta+k+\frac{1}{2},\La+1)$, $(\Delta+k+\frac{1}{2},\La-1)$ and
$(\Delta+k+1,\La)$, respectively. A non-zero $(L_0,J_0)$-weight
vector $v\in M_{\NS^2_+}(\Delta,\La)$ is called a \emph{singular
vector} if $\G_{j} v=0$, for all $j>0$. We call a singular vector
\emph{proper} if it is not a scalar multiple of the highest
weight vector $v_{\Delta,\La}$. Obviously
$M_{\NS^2_+}(\Delta,\La)$ is irreducible if and only if
$M_{\NS^2_+}(\Delta,\La)$ contains no proper singular vector. We
now analyze singular vectors inside $M_{\NS^2_+}(\Delta,\La)$.

\begin{lem}\label{lem1} Let $k\ge 1$ and
suppose that $w=\alpha\pa^kv_{\Delta,\La}
+\beta\pa^{k-1}G^+_{-\frac{1}{2}}G^-_{-\frac{1}{2}}v_{\Delta,\La}$
is a singular vector of $(L_0,J_0)$-weight $(\Delta+k,\La)$ in
$M_{\NS^2_+}(\Delta,\La)$, where $\alpha,\beta\in\C$. Then
$k=1$.  Furthermore any proper singular vector of this form is a
scalar multiple of either
$G^+_{-\frac{1}{2}}G^-_{-\frac{1}{2}}v_{\Delta,\La}$, in which
case $\Delta=-\frac{1}{2}$ and $\La=1$, or
$(-2\pa+G^+_{-\frac{1}{2}}G^-_{-\frac{1}{2}})v_{\Delta,\La}$, in
which case $\Delta=-\frac{1}{2}$ and $\La=-1$.
\end{lem}

\begin{proof}
Note that $w$ is singular if and only if
$J_1w=G^\pm_{\frac{1}{2}}w=0$.  We compute
\begin{align}
G^+_{\frac{1}{2}}w&=(\alpha
k-\beta(2\Delta+\La))\pa^{k-1}G^+_{-\frac{1}{2}}v_{\Delta,\La}=0,\label{generic2}\\
G^-_{\frac{1}{2}}w&=(\alpha
k+\beta(2\Delta-\La+2k))\pa^{k-1}G^-_{-\frac{1}{2}}v_{\Delta,\La}=0,\label{generic3}\\
J_1w&=(\alpha\La k+\beta(2\Delta+\La))\pa^{k-1}v_{\Delta,\La}
+\beta(k-1)\La
\pa^{k-2}G^+_{-\frac{1}{2}}G^-_{-\frac{1}{2}}v_{\Delta,\La}=0.\label{generic4}
\end{align}
But then $\beta\not=0$, since otherwise \eqnref{generic2} would
imply that $k=0$. However, $\beta\not=0$ together with
\eqnref{generic2} and \eqnref{generic3} implies that
\begin{equation}\label{generic5}
2\Delta+k=0.
\end{equation}
Now \eqnref{generic4} gives
\begin{equation}\label{generic6}
\alpha\La k+\beta(2\Delta+\La)=0,\quad \beta(k-1)\La=0.
\end{equation}
Now if $k>1$, then \eqnref{generic6} gives $\La=0$ and
$\Delta=0$. But then $k=0$ by \eqnref{generic2}.  Hence $k=1$ so
that by \eqnref{generic5} we have $\Delta=-\frac{1}{2}$.

Now if $\alpha\not=0$, we have from \eqnref{generic2} and
\eqnref{generic4} $\alpha(1+\La)=0$ and hence $\La=-1$. The first
equation of \eqnref{generic6} then implies that $\alpha+2\beta=0$.

On the other hand if $\alpha=0$, the first equation of
\eqnref{generic6} gives $\La=1$.
\end{proof}

\begin{lem}\label{lem2} Let $k\in\Z_+$.
\begin{itemize}
\item[i.] If $\pa^k G^+_{-\frac{1}{2}}v_{\Delta,\La}$ is a
singular vector of $(L_0,J_0)$-weight $(\Delta+k+\half,\La+1)$,
then $k=0$ and $2\Delta-\La=0$. Furthermore in this case
$G^+_{-\half} v_{\Delta,\La}$ is a singular vector.
\item[ii.] If $\pa^k G^-_{-\frac{1}{2}}v_{\Delta,\La}$ is a
singular vector of $(L_0,J_0)$-weight $(\Delta+k+\half,\La-1)$,
then $k=0$ and $2\Delta+\La=0$. Furthermore in this case
$G^-_{-\half} v_{\Delta,\La}$ is a singular vector.
\end{itemize}
\end{lem}

\begin{proof}
The lemma follows immediately from the following two
equations:
\begin{align*}
G^-_{\frac{1}{2}}\pa^k
G^+_{-\frac{1}{2}}v_{\Delta,\La}&=(2\Delta-\La+2k)\pa^kv_{\Delta,\La}
+k\pa^{k-1}G^+_{-\frac{1}{2}}G^-_{-\frac{1}{2}}v_{\Delta,\La}=0,\\
G^+_{\frac{1}{2}}\pa^k
G^-_{-\frac{1}{2}}v_{\Delta,\La}&=(2\Delta+\La)\pa^kv_{\Delta,\La}
+k\pa^{k-1}G^+_{-\frac{1}{2}}G^-_{-\frac{1}{2}}v_{\Delta,\La}=0.
\end{align*}
\end{proof}

Thus \lemref{lem1} and \lemref{lem2} prove the following.

\begin{prop}\label{singularN=2}
Any proper singular vector in $M_{\NS^2_+}(\Delta,\La)$ is a
scalar multiple of
\begin{itemize}
\item[i.] $G^+_{-\half} v_{\Delta,\La}$, in which case we have $2\Delta-\La=0$.
In the particular case of $\Delta=-\half$ and $\La=-1$ we have in
addition
 $G^-_{-\half} G^+_{-\half} v_{\Delta,\La}$.
\item[ii.] $G^-_{-\half} v_{\Delta,\La}$, in which case we have
$2\Delta+\La=0$. In the particular case of $\Delta=-\half$ and
$\La=1$ we have in addition $G^+_{-\half} G^-_{-\half}
v_{\Delta,\La}$.
\end{itemize}
\end{prop}

Let $N$ be the subspace of $M_{\NS^2_+}(\Delta,\La)$ given by
\begin{align*}
N&=\C[\pa]G^+_{-\frac{1}{2}}v_{\Delta,\La}+\C[\pa]G^-_{-\frac{1}{2}}G^+_{-\frac{1}{2}}v_{\Delta,\La},
\quad {\rm if\ }2\Delta-\La=0{\rm\ and\ }\Lambda\not=0,\\
N&=\C[\pa]G^-_{-\frac{1}{2}}v_{\Delta,\La}+\C[\pa]G^+_{-\frac{1}{2}}G^-_{-\frac{1}{2}}v_{\Delta,\La},
\quad {\rm if\ }2\Delta+\La=0{\rm\ and\ }\Lambda\not=0.
\end{align*}
It follows from \propref{singularN=2} that in either case $N$ is a
submodule of $M_{\NS^2_+}(\Delta,\La)$.

\begin{thm}
The modules $L_{\NS^2_+}(\Delta,\La)$, for $\Delta,\La\in\C$,
form a complete list of non-isomorphic finite (over $\C[\pa]$)
irreducible $K(1,2)_+$-modules.  Furthermore
$L_{\NS^2_+}(\Delta,\La)$ as a $\C[\pa]$-module has rank
\begin{itemize}
\item[i.] $4$, in the case $2\Delta\pm\La\not=0$,
\item[ii.] $2$, in the case $2\Delta\pm\La=0$ and
$2\Delta\mp\La\not=0$,
\item[iii.] $0$, in the case $\Delta=\La=0$.
\end{itemize}
\end{thm}

\begin{proof}
If $2\Delta+\La\not=0$ and $2\Delta-\La\not=0$, then by
\propref{singularN=2} $M_{\NS^2_+}(\Delta,\La)$ contains no
proper singular vector and hence is irreducible.

Suppose that $2\Delta+\La=0$ and $2\Delta-\La\not=0$.  In this
case consider the submodule of $M_{\NS^2_+}(\Delta,\La)$
generated by the singular vector
$G^-_{-\frac{1}{2}}v_{\Delta,\La}$.  This module is precisely $N$
above and hence $M_{\NS^2_+}(\Delta,\La)/N$ is freely generated
over $\C[\pa]$ by $v_{\Delta,\La}$ and $G_{-\half}^+
v_{\Delta,\La}$. We claim that $M_{\NS^2_+}(\Delta,\La)/N$ is
irreducible. The even part of $K(1,2)_+$ is isomorphic to the
semi-direct sum of $\V_+$ (generated by $L_n$) and $\tilde{\G}_+$
(generated by $J_n$), where $\G$ is the one-dimensional Lie
algebra.  We first consider $M_{\NS^2_+}(\Delta,\La)/N$ as a
module over the $\V_+\ltimes\tilde{\G}_+$. The vectors
$v_{\Delta,\La}$ and $G^+_{\half}v_{\Delta,\La}$ have
$(L_0,J_0)$-weights $(\Delta,\La)$ and $(\Delta+\half,\La+1)$,
respectively, and furthermore are both annihilated by $L_n$ and
$J_n$, for $n\ge 1$. Now since $2\Delta+\La=0$ and
$2\Delta-\La\not=0$, we have $(\Delta,\La)\not=(0,0)$ and
$(\Delta+\half,\La+1)\not=(0,0)$. From this it follows that
$M_{\NS^2_+}(\Delta,\La)/N$ as a module over
$\V_+\ltimes\tilde{\G}_+$ is a direct sum of two non-isomorphic
irreducible modules, namely $\C[\pa]v_{\Delta,\La}\cong
L_{\V_+\ltimes\tilde{\G}_+}(\Delta,\La)$ and
$\C[\pa]G^+_{-\frac{1}{2}}v_{\Delta,\La}\cong
L_{\V_+\ltimes\tilde{\G}_+}(\Delta+\frac{1}{2},\La+1)$ (see
\secref{prelim} for notation). But we have
\begin{equation*}
G^-_{\frac{1}{2}}G^+_{-\frac{1}{2}}v_{\Delta,\La}=
(2\Delta-\La)v_{\Delta,\La}\not=0,
\end{equation*}
which implies that as a $K(1,2)_+$-module
$L_{\NS^2_+}(\Delta,\La)$ is irreducible.

The case when $2\Delta-\La=0$ and $2\Delta+\La\not=0$ is
completely analogous and we leave it to the reader.

Finally in the case when $\Delta=\La=0$, both
$G^+_{-\frac{1}{2}}v_{\Delta,\La}$ and
$G^-_{-\frac{1}{2}}v_{\Delta,\La}$ are proper singular vectors.
Now the submodule in $M_{\NS^2_+}(0,0)$ generated by these two
vectors contains
$[G^+_{-\frac{1}{2}},G^-_{-\frac{1}{2}}]v_{\Delta,\La}=2\pa
v_{\Delta,\La}$, and hence has codimension $1$ over $\C$. So the
resulting quotient is the trivial module.
\end{proof}

It follows that every finite irreducible module over the $N=2$
conformal superalgebra is of the form
$L_{\NS^2}(\alpha,\Delta,\La)$, where $\alpha,\Delta,\La\in\C$.
We will write down explicit formulas for the action of the
conformal superalgebra on such irreducible modules in the
generating series form. Since we have already explained in
\secref{prelim} how such formulas can be obtained in general, we
will omit the proofs.

In the case when $2\Delta\pm\La\not=0$ the module
$L_{\NS^2}(\alpha,\Delta,\La)$ is generated freely over $\C[\pa]$
by two even vectors $v,v^{+-}$ and two odd vectors $v^+,v^-$.  We
have the following action on the generators:
\begin{align*}
&L_{\la}v=(\pa+\alpha+\Delta\la)v,\quad
L_{\la}v^\pm=(\pa+\alpha+(\Delta+\half)\la)v^{\pm},\\
&L_{\la}v^{+-}=(\pa+\alpha+(\Delta+1)\la)v^{+-}+(\Delta+\frac{\La}{2})\la^2v,\\
&J_\la v=\La v,\quad J_\la v^{\pm}=(\La\pm 1)v^{\pm},\quad J_\la
v^{+-}=\La v^{+-}+(2\Delta+\La)\la v,\\
&G^{\pm}_{\la}v=v^{\pm},\quad G^{+}_\la v^{+}=G^{-}_\la
v^{-}=0,\quad G^+_{\la}v^-=v^{+-}+(2\Delta+\La)\la v,\\ &G^+_\la
v^{+-}=-\la(2\Delta+\La)v^+,\quad G^-_\la
v^+=(2\pa+2\alpha+\la(2\Delta-\La))v-v^{+-},\\
&G^-_\la v^{+-}=(2\pa+2\alpha+(2\Delta+2-\La)\la)v^-.
\end{align*}

In the case when $2\Delta+\La=0$ but $2\Delta-\La\not=0$ the
module $L_{\NS^2}(\alpha,\Delta,\La)$ is generated freely over
$\C[\pa]$ by one even vector $v$ and one odd vector $v^+$.  The
action is then given by
\begin{align*}
&L_{\la}v=(\pa+\alpha+\Delta\la)v,\quad
L_{\la}v^+=(\pa+\alpha+(\Delta+\half)\la)v^+,\\
&J_\la v=-2\Delta v,\quad J_\la v^+=(-2\Delta + 1)v^+,\quad G^+_{\la}v=v^+,\quad G^+_\la v^+=0,\\
&G^-_\la v=0,\quad G^-_\la
v^+=(2\pa+2\alpha+4\Delta\la) v.\\
\end{align*}

In the case $2\Delta-\La=0$ but $2\Delta+\La\not=0$ the module
$L_{\NS^2}(\alpha,\Delta,\La)$ is generated freely over $\C[\pa]$
by one even vector $v$ and one odd vector $v^-$ with action:
\begin{align*}
&L_{\la}v=(\pa+\alpha+\Delta\la)v,\quad
L_{\la}v^-=(\pa+\alpha+(\Delta+\half)\la)v^-,\\
&J_\la v=2\Delta v,\quad J_\la v^-=(2\Delta- 1)v^-,\quad
G^+_{\la}v=0,\\
&G^+_{\la}v^-=(2\pa+2\alpha+4\Delta \la) v,\quad G^-_\la
v=v^-,\quad G^-_\la v^-=0.
\end{align*}

Finally $L_{\NS^2}(\alpha,0,0)$ is the one-dimensional trivial
module on which $\pa$ acts as the scalar $\alpha$.

\begin{rem} We note that the formulas above are obtained by
first putting $v=v_{\Delta,\La}$,
$v^\pm=G^\pm_{-\half}v_{\Delta,\La}$ and
$v^{+-}=G^{+}_{-\half}G^-_{-\half}v_{\Delta,\La}$ and then
compute the action of the operators $L_n$, $J_m$ and $G^\pm_r$,
for $n\ge -1$, $m\ge 0$ and $r\ge -\half$ on these vector.
Translation into the language of conformal modules is an easy
task using these formulas and we will omit this.  Of course the
parity of the vectors $v,v^\pm,v^{+-}$ in all the examples above
can be reversed. Finally we note that the adjoint module is
isomorphic to $L_{\NS^2}(0,1,0)$.
\end{rem}

\section{Finite irreducible Modules over the $N=3$ conformal superalgebra}\label{N=3}

The $N=3$ superconformal algebra is the formal distribution Lie
superalgebra $K(1,3)$. Letting $\xi_1,\xi_2,\xi_3$ be the three
odd indeterminates $K(1,3)$ is spanned over $\C$ by the following
basis elements ($n\in\Z$ and $r\in\frac{1}{2}+\Z$):
\begin{align*}
&L_n=-\frac{t^{n+1}}{2},\quad H_n=2i\xi_1\xi_2 t^n,\quad
E_n=(-\xi_1\xi_3-i\xi_2\xi_3)t^n,\quad
F_{n}=(\xi_1\xi_3-i\xi_2\xi_3)t^n,\\
&\Psi_r=-\xi_1\xi_2\xi_3 t^{r-\frac{1}{2}}\quad h_r=-2i\xi_3
t^{r+\frac{1}{2}},\quad e_r=(i\xi_1-\xi_2)t^{r+\frac{1}{2}},\quad
f_r=(i\xi_1+\xi_2)t^{r+\frac{1}{2}}.
\end{align*}
Let $\{H,E,F\}$ denote the standard basis of the Lie algebra
$sl_2$ and $\{h,e,f\}$ denote the standard basis of its adjoint
module. Furthermore we let $(\cdot|\cdot)$ denote the
non-degenerate invariant symmetric bilinear form on $sl_2$ with
$(H|H)=2$. Keeping this notation in mind the commutation
relations of $K(1,3)$ are then given as follows (where
$X,Y=H,E,F$ and $x,y=h,e,f$):
\begin{align*}
&[L_m,L_n]=(m-n)L_{m+n},\quad [L_m,X_n]=-nX_{m+n},\quad
[L_m,x_r]=(\frac{m}{2}-r)x_{m+r},\\
&[L_m,\Psi_r]=(-\frac{m}{2}-r)\Psi_{m+r},\quad
[X_m,Y_n]=[X,Y]_{m+n},\quad [X_m,\Psi_r]=0,\\
&[X_m,y_r]=[X,y]_{m+r}+2m(X|Y)\Psi_{m+r},\quad
[x_r,\Psi_s]=-X_{r+s},\quad [\Psi_r,\Psi_s]=0,\\
&[x_r,y_s]=-(r-s)[X,Y]_{r+s}-4(X|Y)L_{r+s},
\end{align*}
where $m,n\in\Z$ and $r,s\in\frac{1}{2}+\Z$. Above we have
written $[X,y]$ for the action of $X$ on $y$.  The eight formal
distributions generating this algebra are given by
$L(z)=\sum_{n\in\Z}L_n z^{-n-2}$, $X(z)=\sum_{n\in\Z}X_n
z^{-n-1}$, $x(z)=\sum_{r\in\frac{1}{2}+\Z}x_r
z^{-r-{\frac{3}{2}}}$ and
$\Psi(z)=\sum_{r\in\frac{1}{2}+\Z}\Psi_r z^{-r-\frac{1}{2}}$. The
corresponding operator product expansions of these fields are
easily derived from \eqnref{generic1}, and so we will omit them.

The annihilation subalgebra $K(1,3)_+$ is equipped with a
$\frac{1}{2}\Z$-gradation of depth $1$, i.e.
$K(1,3)_+=\G=\bigoplus_{j\ge -1}\G_j$, $j\in\frac{1}{2}\Z$, and
its $0$-th graded component $\G_0$ is isomorphic to a copy of
$gl_2\cong sl_2\oplus\C L_0$, with $H_0$, $E_0$ and $F_0$
providing the standard basis for the copy of $sl_2$.

Let $U^{\Delta,\La}$ be the finite-dimensional irreducible
$sl_2$-module of highest weight $\La\in\Z_+$ on which $L_0$ acts
as the scalar $\Delta$.  We let $v_{\Delta,\La}$ be a highest
weight vector in $U^{\Delta,\La}$.  We extend $U^{\Delta,\La}$ to
a module over the subalgebra $\Li_0=\bigoplus_{j\ge 0}\G_j$ in a
trivial way and call this $\Li_0$-module also $U^{\Delta,\La}$.
By \thmref{key} every finite irreducible $\G$-module is a
homomorphic image of $M_{\NS^3_+}(\Delta,\La)={\rm
Ind}_{\Li_0}^{\G}U^{\Delta,\La}$ and furthermore
$M_{\NS^3_+}(\Delta,\La)$ has a unique maximal submodule $N$,
whose irreducible quotient we denote by $L_{\NS^3_+}(\Delta,\La)$.

Note that $M_{\NS^3_+}(\Delta,\La)$ as a module over $sl_2$ is a
direct sum of infinitely many copies of finite-dimensional
irreducible representations.  Since $\pa$ commutes with $E_0$, the
$E_0$-invariants $M_{\NS^3_+}(\Delta,\La)^{E_0}$ is a
$\C[\pa]$-submodule of $M_{\NS^3_+}(\Delta,\La)$, and hence is a
free $\C[\pa]$-module. We can write down explicitly formulas for a
$\C[\pa]$-basis of $M_{\NS^3_+}(\Delta,\La)^{E_0}$.  In the case
when $\La\ge 2$ the following is a $\C[\pa]$-basis:
\begin{align*}
&a_1=v_{\Delta,\La},\quad a_2=e_{-\frac{1}{2}}v_{\Delta,\La},\quad
a_3=(\La h_{-\frac{1}{2}}+2 e_{-\frac{1}{2}}F_0)v_{\Delta,\La},\allowdisplaybreaks\\
&a_4=((\La-1)(\La f_{-\frac{1}{2}}-h_{-\frac{1}{2}}F_0)
-e_{-\frac{1}{2}}F_0^2)v_{\Delta,\La},\allowdisplaybreaks\\
&a_5=e_{-\frac{1}{2}}h_{-\frac{1}{2}}v_{\Delta,\La},\quad a_6=(\La
e_{-\frac{1}{2}}f_{-\frac{1}{2}}
-e_{-\frac{1}{2}}h_{-\frac{1}{2}}F_0)v_{\Delta,\La},\allowdisplaybreaks\\
&a_7=((\La-1)(\La h_{-\frac{1}{2}}f_{-\frac{1}{2}} +4\pa
F_0+2e_{-\frac{1}{2}}f_{-\frac{1}{2}}F_0)
-e_{-\frac{1}{2}}h_{-\frac{1}{2}}F_0^2)v_{\Delta,\La},\allowdisplaybreaks\\
&a_8=(e_{-\frac{1}{2}}h_{-\frac{1}{2}}f_{-\frac{1}{2}}-2\pa
h_{-\frac{1}{2}})v_{\Delta,\La}.
\end{align*}
The cases $\La=0,1$ are similar. Namely, when $\La=1$ we have
$a_4=a_7=0$, and the remaining $6$ vectors form a $\C[\pa]$-basis.
Finally, in the case when $\La=0$, the terms $a_3=a_4=a_6=a_7=0$,
so that $M_{\NS^3_+}(\Delta,0)^{E_0}$ has rank $4$ over
$\C[\pa]$.  (Actually the vectors $a_i$ depend on $\Lambda$, so
it would be more appropriate to write something like $a_i^\La$
instead of just $a_i$. However, from the context it will always
be clear what $\La$ is, so that it is safe to adopt the simpler
notation of $a_i$.) We denote the coefficient of $v_{\Delta,\La}$
in the expression $a_i$ by $u^\La_i$ so that we have
$a_i=u^\La_iv_{\Delta,\La}$, for $i=1,\ldots,8$. For example
$u^\La_1=1$, while $u^\La_2=e_{-\frac{1}{2}}$ etc.  We note that
finding all vectors in $M_{\NS^3_+}(\Delta,\La)^{E_0}$ above
amounts essentially to decomposing tensor products of irreducible
representations of $sl_2$ and then finding the corresponding
highest weight vectors of the irreducible components.

Similarly we call a non-zero $(L_0,H_0)$-weight vector $v$ in
$M_{\NS^3_+}(\Delta,\La)$ \emph{singular} if $v\in
M_{\NS^3_+}(\Delta,\La)^{E_0}$ and $\G_{j}v=0$, for all $j>0$. As
before a singular vector is called \emph{proper} if it is not a
scalar multiple of $v_{\Delta,\La}$. Evidently
$M_{\NS^3_+}(\Delta,\La)$ is irreducible if and only if
$M_{\NS^3_+}(\Delta,\La)$ contains no proper singular vector. Our
first objective is to classify singular vectors inside
$M_{\NS^3_+}(\Delta,\La)$.

\begin{prop}\label{singularN=3}
Any proper singular vector in $M_{\NS^3_+}(\Delta,\La)$ is of the
form ($\alpha\in\C$ with $\alpha\not=0$)
\begin{itemize}
\item[i.] $\alpha a_2$, if $4\Delta-\La=0$,
\item[ii.] $\alpha a_4$, if $4\Delta+\La+2=0$ and $\La\ge 2$,
\item[iii.] $\alpha a_6$, if $4\Delta+\La+2=0$ and $\La=1$.
\end{itemize}
\end{prop}

\begin{rem}\label{N=3rem1}
The proof of the proposition is a straightforward, albeit a
tedious, calculation.  We will not give the details here, but
instead just point out that a weight vector $v\in
M_{\NS^3_+}(\Delta,\La)^{E_0}$ is singular if and only if
$f_{\frac{1}{2}}$ and $\Psi_{\frac{1}{2}}$ annihilates $v$.  This
fact simplifies the calculation significantly.
\end{rem}

From \propref{singularN=3} one obtains immediately the following.

\begin{cor}
Suppose that $(\Delta,\La)$ does not satisfy either
$4\Delta-\La=0$ or $4\Delta+\La+2=0$ and $\La\ge 1$. Then
$L_{\NS^3_+}(\Delta,\La)=M_{\NS^3_+}(\Delta,\La)$ is an
irreducible $K(1,3)_+$-module of rank $8\La+8$ over $\C[\pa]$.
\end{cor}

\begin{prop}\label{N=3case1}
Suppose that $4\Delta-\La=0$. Then $L_{\NS^3_+}(\Delta,\La)$ is a
free $\C[\pa]$-module of rank $4\La$.
\end{prop}

\begin{proof}
By \propref{singularN=3} $a_2$ is a singular vector in
$M_{\NS^3_+}(\frac{\La}{4},\La)$ of $H_0$-weight $\La+2$. Consider
$N$, the $\G$-submodule generated by $a_2$. Then we have
$N=U(\G_{-})V_2$, where $V_2$ is the irreducible $sl_2$-submodule
generated by $a_2$. Note that the map
$v_{\frac{\La}{4}+\frac{1}{2},\La+2}\rightarrow a_2$ extends
uniquely to an epimorphism of $K(1,3)_+$-modules from
$M_{\NS^3_+}(\frac{\La}{4}+\frac{1}{2},\La+2)$ to $N$.  In
particular it is an $sl_2$-module epimorphism.  Now both modules
are completely reducible $sl_2$-modules and hence this map sends
$E_0$-invariants onto $E_0$-invariants. Since
$M_{\NS^3_+}(\frac{\La}{4}+\frac{1}{2},\La+2)^{E_0}$ is generated
over $\C[\pa]$ by
$\{u^{\La+2}_iv_{\frac{\La}{4}+\frac{1}{2},\La+2}|1\le i\le 8\}$,
it follows that $N^{E_0}$ is generated over $\C[\pa]$ by
$\{u^{\La+2}_ia_2|1\le i\le 8\}$. Now $N^{E_0}$ is a
$\C[\pa]$-submodule of $M_{\NS^3_+}(\frac{\La}{4},\La)$, since
$[\pa,E_0]=0$.  Thus it is a free $\C[\pa]$-submodule generated
by $\{u^{\La+2}_ia_2|1\le i\le 8\}$. We compute
\begin{align*}
&u^{\La+2}_1a_2=a_2,\quad u^{\La+2}_2a_2=0,\quad u^{\La+2}_3a_2=-(\La+4)a_5,\\
&u^{\La+2}_4a_2=-(\La+3)a_6-4(\La+1)(\La+3)\pa a_1,\quad u^{\La+2}_5 a_2=0,\\
&u^{\La+2}_6a_2=-4(\La+3)\pa a_2,\quad
u^{\La+2}_7a_2=(\La+3)(\La+2)a_8+2(\La+3)\pa a_3,\\
&u^{\La+2}_8a_2=-2\pa a_5.
\end{align*}

By inspection it is clear that the following is a set of
$\C[\pa]$-generators for $N^{E_0}$.
\begin{equation*}
S^\La=\{a_2,a_5,a_6+4(\La+1)\pa a_1,a_8+\frac{2}{\La+2}\pa a_3\}.
\end{equation*}
First consider the case when $\La\ge 2$. It follows from the
description of $S^\la$ above that $\{a_1,a_3,a_4,a_7\}$ is a
$\C[\pa]$-basis for the $E_0$-invariants of the quotient
$M_{\NS^3_+}(\frac{\La}{4},\La)/N$. Since $a_1$ and $a_3$ both
have $H_0$-weight $\La$, they generate two copies of the
irreducible $sl_2$-module of dimension $\La+1$. On the other hand
$a_4$ and $a_7$ both have weight $\La-2$, and so they generate
two copies of the irreducible $sl_2$-module of dimension $\La-1$.
Thus $M_{\NS^3_+}(\frac{\La}{4},\La)/N$ is a free
$\C[\pa]$-module of rank $2(\La+1)+2(\La-1)=4\La$.  So in order
to complete the proof it remains to show that
$M_{\NS^3_+}(\frac{\La}{4},\La)/N$ is irreducible.

Note that $L_n$, $n\ge -1$, together with $E_0,H_0,F_0$ generate a
copy of $\V_+\oplus sl_2$ and so we may consider
$M_{\NS^3_+}(\frac{\La}{4},\La)/N$ as a module over $\V_+\oplus
sl_2$. By parity consideration $M_{\NS^3_+}(\frac{\La}{4},\La)/N$
is a direct sum of two $(\V_+\oplus sl_2)$-modules, namely
$(M_{\NS^3_+}(\frac{\La}{4},\La)/N)_{\bar{0}}=\C[\pa]V_1+\C[\pa]V_7$
and
$(M_{\NS^3_+}(\frac{\La}{4},\La)/N)_{\bar{1}}=\C[\pa]V_3+\C[\pa]V_4$,
where $V_i$ is the irreducible $sl_2$-module generated by $a_i$.
It is subject to a direct verification that $L_n$, for $n\ge 1$,
annihilates the vectors $a_1,a_3,a_4,a_7$ (in fact one only needs
to check that $L_1 a_7=0$, others being trivial) and hence
$M_{\NS^3_+}(\frac{\La}{4},\La)/N$ as a $\V_+\oplus sl_2$-module
is a direct sum the following four non-isomorphic irreducible
modules: $\C[\pa]V_3\cong
L_{\V_+}(\frac{\La}{4}+\frac{1}{2})\boxtimes U^\La$,
$\C[\pa]V_4\cong L_{\V_+}(\frac{\La}{4}+\frac{1}{2})\boxtimes
U^{\La-2}$, $\C[\pa]V_1\cong L_{\V_+}(\frac{\La}{4})\boxtimes
U^\La$ and $\C[\pa]V_7\cong L_{\V_+}(\frac{\La}{4}+1)\boxtimes
U^{\La-2}$, where we denote by $U^\mu$ the irreducible
$sl_2$-module of highest weight $\mu$.  Now we compute
\begin{align}
&\Psi_{\frac{1}{2}}a_3=-\La(\La+2)a_1,\quad
f_{\frac{1}{2}} a_4=(2\La+2)F_0^2a_1\label{N=3aux1}\\
&E_{1}a_7=2\La(\La-1)(2\La+2)a_1,\nonumber
\end{align}
from which it follows that we may go from each irreducible
$\V_+\oplus sl_2$-component of $M_{\NS^3_+}(\frac{\La}{4},\La)/N$
to the irreducible component containing the highest weight
vectors, and hence the module $M_{\NS^3_+}(\frac{\La}{4},\La)/N$
is irreducible.

Now if $\La=1$ the vectors $a_4$ and $a_7$ are both zero.
Therefore the quotient
$M_{\NS^3_+}(\frac{\La}{4},\La)/N=\C[\pa]V_1\oplus\C[\pa]V_3$.
But then the first identity in \eqnref{N=3aux1} shows that
$M_{\NS^3_+}(\frac{\La}{4},\La)/N$ is irreducible.  The rank of
$L_{\NS^3_+}(\frac{\La}{4},\La)$ is then $2(\La+1)=4\La$ in the
case $\La=1$.

Finally when $\La=0$, the vectors $a_3,a_4,a_6,a_7=0$, so that
$S^\La$ reduces to $\{a_2,a_5,\pa a_1,a_8\}$. Therefore
$M_{\NS^3_+}(0,0)/N=\C a_1$ is the trivial module and so has rank
$0$.
\end{proof}

\begin{prop}\label{N=3case2}
Suppose that $4\Delta+\La+2=0$ and $\La\ge 1$.  Then
$L_{\NS^3_+}(\Delta,\La)$ is a free $\C[\pa]$-module of rank
$4\La+8$.
\end{prop}

\begin{proof}
By \propref{singularN=3} $a_4$ is a singular vector of
$M_{\NS^3_+}(-\frac{\La+2}{4},\La)$ of $H_0$-weight $\La-2$. Let
$N$ denote the $\G$-submodule generated by $a_4$.

Consider first the case $\La\ge 4$. As in the proof of
\propref{singularN=3} $N^{E_0}$ is a free $\C[\pa]$-module
generated over $\C[\pa]$ by $\{u^{\La-2}_ia_4|1\le i\le 8\}$.  We
compute
\begin{align*}
&u^{\La-2}_1a_4=a_4,\quad u^{\La-2}_2a_4=(\La-1)a_6,\quad
u^{\La-2}_3a_4=(\La-2) a_{7},\quad u^{\La-2}_4a_4=0,\\
&u^{\La-2}_5a_4=\La(\La-1)a_8+2(\La-1)\pa a_3,\quad
u^{\La-2}_6a_4=0,\quad u^{\La-2}_7a_4=0,\\ &u^{\La-2}_8a_4=-2\pa
a_7.
\end{align*}
This implies that the set
$S^\La=\{a_4,(\La-1)a_6,(\La-2)a_7,\La(\La-1)a_8+2(\La-1)\pa
a_3,\pa a_7 \}$ generates $N^{E_0}$ over $\C[\pa]$ and so
$\{a_1,a_2,a_3,a_5\}$ is a $\C[\pa]$-basis for
$(M_{\NS^3_+}(-\frac{\La+2}{4},\La)/N)^{E_0}$ in the case when
$\La\ge 4$.

Next consider the case $\La=3$.  In this case, letting $N$ be as
before, $N^{E_0}$ is generated over $\C[\pa]$ by
$\{u^{\La-2}_ia_4|1\le i\le 8,i\not=4,7\}$. Hence it follows from
the above formulas that again $\{a_1,a_2,a_3,a_5\}$ is a
$\C[\pa]$-basis for $(M_{\NS^3_+}(-\frac{\La+2}{4},\La)/N)^{E_0}$.

In the case when $\La=2$ we let $N'$ denote the module generated
by $a_4$. It follows that the vectors
$\{u^{\La-2}_1a_4,u^{\La-2}_2a_4,u^{\La-2}_5a_4,u^{\La-2}_8a_4\}$
generate ${N'}^{E_0}$ over $\C[\pa]$ so that $S^\La=\{a_4,a_6,
a_8+\pa a_3,\pa a_7\}$ generate $N^{E_0}$. Hence
$(M_{\NS^3_+}(-\frac{\La+2}{4},\La)/N')^{E_0}$ contains in
addition a one-dimensional (over $\C$) subspace spanned by $a_7$.
However, $\pa a_7=0$ in $(M_{\NS^3_+}(-\frac{\La+2}{4},\La)/N')$
and hence it is a $\G$-invariant by \remref{invariant}.  In this
case we set $N=N'+\C a_7$ so that the quotient module
$(M_{\NS^3_+}(-\frac{\La+2}{4},\La)/N)^{E_0}$ is again generated
over $\C[\pa]$ by $\{a_1,a_2,a_3,a_5\}$.

Now $a_1$ and $a_3$ have $H_0$-weight $\La$, while $a_2$ and
$a_5$ have $H_0$-weight $\La+2$.  Thus
$M_{\NS^3_+}(-\frac{\La+2}{4},\La)/N$ has rank
$2(\La+1)+2(\La+3)=4\La+8$ over $\C[\pa]$.  So it remains to show
that $M_{\NS^3_+}(-\frac{\La+2}{4},\La)/N$ is irreducible.

Again we consider $M_{\NS^3_+}(-\frac{\La+2}{4},\La)/N$ as a
module over $\V_+\oplus sl_2$.  It is easy to check that $L_n$,
$n\ge 1$, annihilates $a_1,a_2,a_3,a_5$.  (Again one really only
needs to check that $L_1a_5=0$.)  Thus it follows that in the
case of $\La\ge 3$ that $M_{\NS^3_+}(-\frac{\La+2}{4},\La)/N$ is a
direct sum of the following four non-isomorphic irreducible
$\V_+\oplus sl_2$-modules: $\C[\pa]V_1\cong
L_{\V_+}(-\frac{\La+2}{4})\boxtimes U^\La$, $\C[\pa]V_2\cong
L_{\V_+}(-\frac{\La}{4})\boxtimes U^{\La+2}$, $\C[\pa]V_3\cong
L_{\V_+}(-\frac{\La}{4})\boxtimes U^\La$ and $\C[\pa]V_5\cong
L_{\V_+}(-\frac{\La-2}{4})\boxtimes U^{\La+2}$, where as before
$U^\mu$ stands for the irreducible $sl_2$-module of highest
weight $\mu$ and $V_i$ is the $sl_2$-submodule generated by the
vector $a_i$. Now we compute
\begin{equation}\label{N=3generic11}
f_{\frac{1}{2}}a_2=2(\La+1)a_1,\quad
\Psi_{\frac{1}{2}}a_3=-\La(\La+2)a_1,\quad F_1a_5=-4(\La+1)a_1,
\end{equation}
from which again it follows that we may go from any irreducible
$\V_+\oplus sl_2$-component of
$M_{\NS^3_+}(-\frac{\La+2}{4},\La)/N$  to the component containing
the highest weight vectors, and hence
$M_{\NS^3_+}(-\frac{\La+2}{4},\La)/N$ is irreducible.

As for the case $\La=2$ we have
$M_{\NS^3_+}(-\frac{\La+2}{4},\La)/N$ as a $\V_+\oplus
sl_2$-module is also a direct sum of the
$\C[\pa]V_1\oplus\C[\pa]V_2\oplus\C[\pa]V_3\oplus\C[\pa]V_5$. The
first three modules, as in the case of $\La\ge 3$ are irreducible.
However, $\C[\pa]V_5$ contains a unique irreducible submodule
generated by the vector $\pa a_5$, which is isomorphic to
$L_{\V_+}(-1)\boxtimes U^{\La+2}$.  But then
\eqnref{N=3generic11} together with the fact that $$F_2\pa a_1=-24
a_1$$ shows that $M_{\NS^3_+}(-\frac{\La+2}{4},\La)/N$ is
irreducible in this case as well.

In the case when $\La=1$ we have by \propref{singularN=3} that
$a_6$ is the unique (up to a scalar) singular vector inside
$M_{\NS^3_+}(-\frac{3}{4},1)$.  Let $N$ denote the $\G$-submodule
generated by $a_6$.  Since $a_6$ has $H_0$-weight $1$, $N^{E_0}$
is the free $\C[\pa]$-module generated by $\{u^{\La}_ia_6|1\le
i\le 8,i\not=4,7\}$.  We have
\begin{align*}
&u^{\La}_1a_6=a_6,\quad u^{\La}_2a_6=0,\quad u^{\La}_3a_6=-3a_8-6\pa a_3,\\
&u^{\La}_5a_6=0,\quad u^{\La}_6a_6=-8\pa a_6,\quad
u^{\La}_8a_6=-2\pa a_8-4\pa^2 a_3,
\end{align*}
from which it follows that $N^{E_0}$ is generated over $\C[\pa]$
by $S^\La=\{a_6,a_8+2\pa a_3\}$.  Since $a_4=a_7=0$ in this
situation, we see that $(M_{\NS^3_+}(-\frac{3}{4},1)/N)^{E_0}$ is
generated over $\C[\pa]$ by the vectors $a_1,a_2,a_3,a_5$, just
as in the case $\La\ge 2$.  Now the exact same argument as in the
$\La\ge 2$ case shows that $M_{\NS^3_+}(-\frac{3}{4},1)/N$ is
irreducible and has rank $4\La+8$ over $\C[\pa]$.
\end{proof}

We summarize the work in this section in the following theorem.

\begin{thm}
The modules $L_{\NS^3_+}(\Delta,\La)$, for $\Delta\in\C$ and
$\La\in\Z_+$, form a complete list of non-isomorphic finite (over
$\C[\pa]$) irreducible $K(1,3)_+$-modules.  Furthermore
$L_{\NS^3_+}(\Delta,\La)$ as a $\C[\pa]$-module has rank
\begin{itemize}
\item[i.] $4\La$, in the case $4\Delta-\La=0$,
\item[ii.] $4\La+8$, the case $4\Delta+\La+2=0$
and $\La\ge 1$.
\item[iii.] $8\La+8$, in all other cases.
\end{itemize}
Furthermore the $\C[\pa]$-rank of
$L_{\NS^3_+}(\Delta,\La)_{\bar{0}}$ equals the $\C[\pa]$-rank of
$L_{\NS^3_+}(\Delta,\La)_{\bar{1}}$ in all cases.
\end{thm}

\begin{rem}
Translating the above theorem back into the languages of modules
over conformal superalgebras and of conformal modules is now a
straightforward task.   We thus have proved that all finite
irreducible modules over the $N=3$ conformal superalgebra are of
the form $L_{\NS^3}(\alpha,\Delta,\La)$, where
$\alpha,\Delta\in\C$ and $\La\in\Z_+$.  The definition of these
modules and also the action of the $N=3$ conformal superalgebra
on them are quite easy to obtain from our explicit description of
a $\C[\pa]$-basis of these modules.  To do so would however take
up quite a significant portion of space, and thus we leave this
task to the interested reader.  We only remark that the adjoint
module is isomorphic to $L_{\NS^3}(0,\half,0)$.
\end{rem}

\section{Finite irreducible Modules over the ``small'' $N=4$ conformal superalgebra}\label{N=4}

The ``small'' $N=4$ superconformal algebra is the following
subalgebra of $K(1,4)$:  Let $\xi_1,\xi_2,\xi_3,\xi_4$ denote
four odd indeterminates generating the Grassmann superalgebra
$\Lambda(4)$. For a monomial $\xi_I$ in $\Lambda(4)$ we let
$\xi_I^*$ be its Hodge dual, i.e.~the unique monomial in
$\Lambda(4)$ such that $\xi_I\xi_I^*=\xi_1\xi_2\xi_3\xi_4$.  Then
the small $N=4$ superconformal algebra is isomorphic to any of
the following two subalgebras in $K(1,4)$ spanned by the
following basis elements ($n\in\Z$, $r\in\frac{1}{2}+\Z$,
$\beta^2=1$) \cite{CK2}:
\begin{align*}
&L^{\beta}_n=-\frac{1}{2}(t^{n+1}+{\beta}n(n+1)\xi_1\xi_2\xi_3\xi_4
t^{n-1}),\\
&H^{\beta}_n=i(\xi_1\xi_2-{\beta}\xi_3\xi_4)t^n,\\
&E^{\beta}_n=\frac{1}{2}(-\xi_1\xi_3-{\beta}\xi_2\xi_4-i\xi_2\xi_3+i{\beta}\xi_1\xi_4)t^n,\\
&F^{\beta}_n=\frac{1}{2}(\xi_1\xi_3+{\beta}\xi_2\xi_4-i\xi_2\xi_3+i{\beta}\xi_1\xi_4)t^n,\allowdisplaybreaks\\
&G^{-+{\beta}}_r=\frac{1}{\sqrt{2}}((\xi_3+i\xi_4)t^{r+\frac{1}{2}}-
{\beta}(r+\frac{1}{2})(\xi^*_3+i\xi^*_4)t^{r-\frac{1}{2}}),\allowdisplaybreaks\\
&G^{++{\beta}}_r=\frac{1}{\sqrt{2}}((\xi_1+i\xi_2)t^{r+\frac{1}{2}}-
{\beta}(r+\frac{1}{2})(\xi^*_1+i\xi^*_2)t^{r-\frac{1}{2}}),\allowdisplaybreaks\\
&G^{+-{\beta}}_r=\frac{1}{\sqrt{2}}((\xi_3-i\xi_4)t^{r+\frac{1}{2}}-
{\beta}(r+\frac{1}{2})(\xi^*_3-i\xi^*_4)t^{r-\frac{1}{2}}),\allowdisplaybreaks\\
&G^{--{\beta}}_r=\frac{1}{\sqrt{2}}((i\xi_2-\xi_1)t^{r+\frac{1}{2}}-
{\beta}(r+\frac{1}{2})(i\xi^*_2-\xi^*_1)t^{r-\frac{1}{2}}).
\end{align*}

As before let $\{H,E,F\}$ denote the standard basis of the Lie
algebra $sl_2$ and $\{G^{++},G^{-+}\}$ denote the standard basis
of its standard module, i.e. $H\cdot G^{++}=G^{++}$, $H\cdot
G^{-+}=-G^{-+}$, $E\cdot G^{++}=F\cdot G^{-+}=0$, $F\cdot
G^{++}=G^{-+}$ and $E\cdot G^{-+}=G^{++}$. Likewise
$\{G^{+-},G^{--}\}$ also denotes a copy of the standard basis of
the standard $sl_2$-module with actions $H\cdot G^{+-}=G^{+-}$,
$H\cdot G^{--}=-G^{--}$, $E\cdot G^{+-}=F\cdot G^{--}=0$, $F\cdot
G^{+-}=G^{--}$ and $E\cdot G^{--}=G^{+-}$. With this notation in
mind the commutation relations are then given as follows (where
$X,Y=H,E,F$ and $x,y=G^{++},G^{-+},G^{+-},G^{--}$):
\begin{align*}
&[L^{\beta}_m,L^{\beta}_n]=(m-n)L^{\beta}_{m+n},\quad
[L^{\beta}_m,X^{\beta}_n]=-nX^{\beta}_{m+n},\quad
[L^{\beta}_m,x^{\beta}_r]=(\frac{m}{2}-r)x^{\beta}_{m+r},\\
&[X^{\beta}_m,Y^{\beta}_n]=[X,Y]^{\beta}_{m+n},\quad [X^{\beta}_m,y^{\beta}_r]=(X\cdot y)^{\beta}_{m+r},
\quad [x^{\beta}_r,x^{\beta}_s]=0,\allowdisplaybreaks\\
&[G^{++\beta}_r,G^{+-\beta}_s]=(r-s)(1+\beta)E^\beta_{r+s},\quad
[G^{++\beta}_r,G^{-+\beta}_s]=(r-s)(1-\beta)E^\beta_{r+s},\allowdisplaybreaks\\
&[G^{++\beta}_r,G^{--\beta}_s]=-(r-s)H^\beta_{r+s}-2L^\beta_{r+s},\
[G^{+-\beta}_r,G^{-+\beta}_s]=(r-s)\beta
H^\beta_{r+s}+2L^\beta_{r+s},\allowdisplaybreaks\\
&[G^{+-\beta}_r,G^{--\beta}_s]=-(r-s)(1-\beta)F^\beta_{r+s},\quad
[G^{-+\beta}_r,G^{--\beta}_s]=-(r-s)(1+\beta)F^\beta_{r+s},
\end{align*}
where $m,n\in\Z$ and $r,s\in\frac{1}{2}+\Z$. The eight formal
distributions generating this algebra are given by
$L^\beta(z)=\sum_{n\in\Z}L^{\beta}_n z^{-n-2}$,
$X^\beta(z)=\sum_{n\in\Z}X^{\beta}_n z^{-n-1}$,
$x^\beta(z)=\sum_{r\in\frac{1}{2}+\Z}x^{\beta}_r
z^{-r-{\frac{3}{2}}}$. The operator product expansions of these
fields are easily derived using \eqnref{generic1}.

We will denote the ``small'' $N=4$ superconformal algebra simply
by $SK(1,4)$ and assume for the rest of this section that we have
chosen its realization as the subalgebra of $K(1,4)$ with
$\beta=1$ for future computational purposes. For simplicity we
will drop the superscript $\beta$ and write $L_n$ for
$L^{\beta}_n$ etc.~when we mean $\beta=1$.

The annihilation subalgebra $\G=SK(1,4)_+$ of $SK(1,4)$ is
equipped with a $\frac{1}{2}\Z$-gradation of depth $1$, i.e.
$\G=\bigoplus_{j\ge -1}\G_j$, $j\in\frac{1}{2}\Z$, and its $0$-th
graded component $\G_0$ is isomorphic to a copy of $gl_2\cong
sl_2\oplus\C L_0$, with $H_0$, $E_0$ and $F_0$ providing the
standard basis of the copy of $sl_2$. Again we let
$U^{\Delta,\La}$ be the finite-dimensional irreducible
$sl_2$-module of highest weight $\La\in\Z_+$ on which $L_0$ acts
as the scalar $\Delta$ and $v_{\Delta,\La}$ be a highest weight
vector in $U^{\Delta,\La}$. As in the case of $K(1,3)_+$, we may
extend $U^{\Delta,\La}$ to a module over the subalgebra
$\Li_0=\bigoplus_{j\ge 0}\G_j$ trivially and call this
$\Li_0$-module also $U^{\Delta,\La}$. Again \thmref{key} tells us
that every finite irreducible $\G$-module is the quotient of
$M_{\NS^4_+}(\Delta,\La)={\rm Ind}_{\Li_0}^{\G}U^{\Delta,\La}$ by
its unique maximal submodule, for some $\Delta\in\C$ and
$\La\in\Z_+$. We denote the unique irreducible quotient by
$L_{\NS^4_+}(\Delta,\La)$ so that every finite irreducible
$SK(1,4)_+$-module is of the form $L_{\NS^4_+}(\Delta,\La)$, for
$\Delta\in\C$ and $\La\in\Z_+$.

Now $M_{\NS^4_+}(\Delta,\La)$ is completely reducible as a module
over $sl_2=\C H_0+\C E_0+\C F_0$, and the subspace of
$E_0$-invariants $M_{\NS^4_+}(\Delta,\La)^{E_0}$ is a free
$\C[\pa]$-submodule of $M_{\NS^4_+}(\Delta,\La)$ due to
$[E_0,\pa]=0$. We write down explicit formulas for a
$\C[\pa]$-basis of $M_{\NS^4_+}(\Delta,\La)^{E_0}$, which in the
case when $\La\ge 2$ takes the following form:
\begin{align*}
&a_1=v_{\Delta,\La},\quad
a_2=G^{++}_{-\frac{1}{2}}v_{\Delta,\La},\quad
a_3=G^{+-}_{-\frac{1}{2}}v_{\Delta,\La},\quad
a_4=(\La G^{-+}_{-\half}-G^{++}_{-\half} F_0)v_{\Delta,\La},\allowdisplaybreaks\\
&a_5=(-\La
G^{--}_{-\frac{1}{2}}+G^{+-}_{-\frac{1}{2}}F_0)v_{\Delta,\La},\quad
a_6=G^{-+}_{-\frac{1}{2}}G^{++}_{-\frac{1}{2}}v_{\Delta,\La},\quad
a_7=G^{+-}_{-\frac{1}{2}}G^{--}_{-\frac{1}{2}}v_{\Delta,\La},\allowdisplaybreaks\\
&a_8=G^{++}_{-\frac{1}{2}}G^{+-}_{-\frac{1}{2}}v_{\Delta,\La},\quad
a_9=(G^{-+}_{-\frac{1}{2}}G^{+-}_{-\frac{1}{2}}-
G^{++}_{-\frac{1}{2}} G^{--}_{-\frac{1}{2}})v_{\Delta,\La},\allowdisplaybreaks\\
&a_{10}=(-\La G^{++}_{-\frac{1}{2}}G^{--}_{-\frac{1}{2}}+
G^{++}_{-\frac{1}{2}}G^{+-}_{-\frac{1}{2}} F_0)v_{\Delta,\La},\allowdisplaybreaks\\
&a_{11}=((\La-1)(-\La G^{-+}_{-\frac{1}{2}}G^{--}_{-\frac{1}{2}}+
G^{-+}_{-\frac{1}{2}}G^{+-}_{-\frac{1}{2}}
F_0+G^{++}_{-\frac{1}{2}}G^{--}_{-\frac{1}{2}} F_0)
-G^{++}_{-\frac{1}{2}} G^{+-}_{-\frac{1}{2}}F_0^2)v_{\Delta,\La},\allowdisplaybreaks\\
&a_{12}=G^{-+}_{-\frac{1}{2}}G^{++}_{-\frac{1}{2}}G^{+-}_{-\frac{1}{2}}
v_{\Delta,\La},\quad
a_{13}=G^{++}_{-\frac{1}{2}}G^{+-}_{-\frac{1}{2}}G^{--}_{-\frac{1}{2}}
v_{\Delta,\La},\allowdisplaybreaks\\
&a_{14}=G^{-+}_{-\frac{1}{2}}G^{++}_{-\frac{1}{2}}(-\La
G^{--}_{-\frac{1}{2}}+
G^{+-}_{-\frac{1}{2}}F_0)v_{\Delta,\La},\allowdisplaybreaks\\
&a_{15}=(-\La
G^{-+}_{-\frac{1}{2}}G^{+-}_{-\frac{1}{2}}G^{--}_{-\frac{1}{2}}
+G^{++}_{-\frac{1}{2}}G^{+-}_{-\frac{1}{2}}G^{--}_{-\frac{1}{2}}F_0
)v_{\Delta,\La},\allowdisplaybreaks\\
&a_{16}=G^{-+}_{-\frac{1}{2}}G^{++}_{-\frac{1}{2}}
G^{+-}_{-\frac{1}{2}}G^{--}_{-\frac{1}{2}}v_{\Delta,\La}.
\end{align*}
Now in the case when $\La=1$ we have $a_{11}=0$ so that the
remaining $15$ vectors form a $\C[\pa]$-basis for
$M_{\NS^4_+}(\Delta,\La)^{E_0}$, while in the case when $\La=0$
we have $a_4=a_5=a_{10}=a_{11}=a_{14}=a_{15}=0$, so that
$M_{\NS^4_+}(\Delta,0)^{E_0}$ has rank $10$ over $\C[\pa]$.  As in
\secref{N=3} we denote the coefficient of $v_{\Delta,\La}$ in the
expression $a_i$ by $u^\La_i$ so that we have
$a_i=u^\La_iv_{\Delta,\La}$, for $i=1,\ldots,16$.

\emph{Singular vectors} are then defined to be non-zero
$(L_0,H_0)$-weight vectors $v\in M_{\NS^4_+}(\Delta,\La)^{E_0}$
with $\G_{j}v=0$, for all $j>0$. Similarly we define \emph{proper
singular vectors}. Our approach is analogous to the one of
\secref{N=3}, that is first to analyze singular vectors inside
$M_{\NS^4_+}(\Delta,\La)$.  This is given by the following
proposition, whose proof is again a straightforward calculation,
which admittedly is rather tedious.

\begin{prop}\label{singularN=4}
A complete list of proper singular vectors inside
$M_{\NS^4_+}(\Delta,\La)$ is given by:
\begin{itemize}
\item[i.] $2\Delta-\La=0$.
\begin{itemize}
\item[a.] $\alpha a_2+\beta a_3$, $(\alpha,\beta)\not=(0,0)$,
\item[b.] $\alpha a_8$, $\alpha\not=0$.
\end{itemize}
\item[ii.] $2\Delta+\La+2=0$ and $\La\ge 2$.
\begin{itemize}
\item[a.] $\alpha a_4+\beta a_5$, $(\alpha,\beta)\not=(0,0)$,
\item[b.] $\alpha a_{11}$, $\alpha\not=0$.
\end{itemize}
\item[iii.] $2\Delta+\La+2=0$ and $\La=1$.
\begin{itemize}
\item[a.] $\alpha a_4+\beta a_5$, $(\alpha,\beta)\not=(0,0)$,
\item[b.] $\alpha a_{14}+\beta(a_{15}-2\pa a_5)$, $(\alpha,\beta)\not=(0,0)$,
\item[c.] $\alpha (a_{16}-2\pa a_{10})$, $\alpha\not=0$.
\end{itemize}
\item[iv.] $2\Delta+\La+2=0$ and $\La=0$.
\begin{itemize}
\item[a.] $\alpha a_6+\beta a_7+\gamma(a_9-2\pa a_1)$, $(\alpha,\beta,\gamma)\not=(0,0,0)$,
\item[b.] $\alpha a_{13}+\beta(a_{12}+2\pa a_2)$,
$(\alpha,\beta)\not=(0,0)$.
\end{itemize}
\end{itemize}
\end{prop}

\begin{rem}\label{N=4rem1}
We note that in order to check that a weight vector $v\in
M_{\NS^4_+}(\Delta,\La)^{E_0}$ is singular, it is enough to check
that $v$ is annihilated by $F_1$, $G^{-+}_{\frac{1}{2}}$ and
$G^{--}_{\frac{1}{2}}$.
\end{rem}

\begin{cor}\label{N=4generic}
Suppose that $(\Delta,\La)$ does not satisfy either
$2\Delta-\La=0$ or $2\Delta+\La+2=0$. Then
$L_{\NS^4_+}(\Delta,\La)=M_{\NS^4_+}(\Delta,\La)$ is an
irreducible $SK(1,4)_+$-module of rank $16\La+16$ over $\C[\pa]$.
\end{cor}

\begin{prop}
Suppose that $2\Delta-\La=0$. Then $L_{\NS^4_+}(\Delta,\La)$ is a
free $\C[\pa]$-module of rank $4\La$.
\end{prop}

\begin{proof}
By \propref{singularN=4} $a_2$ and $a_3$ are singular vectors in
$M_{\NS^4_+}(\frac{\La}{2},\La)$. Consider $N_2$ and $N_3$, the
$\G$-submodules generated by $a_2$ and $a_3$, respectively, and
let $N=N_2+N_3$.  Then we have $N_2=U(\G_{-})V_2$ and
$N_3=U(\G_{-})V_3$, where $V_2$ and $V_3$ are the irreducible
$sl_2$-submodules generated by $a_2$ and $a_3$, respectively.
Let's first compute $N^{E_0}_2$. Since the $H_0$-weight of $a_2$
is $\La+1$, we know that $N_2^{E_0}$ is a free $\C[\pa]$-module
generated over $\C[\pa]$ by $\{u^{\La+1}_ia_2|1\le i\le 16\}$. We
have
\begin{align*}
&u^{\La+1}_1a_2=a_2,\quad u^{\La+1}_2a_2=0,\quad
u^{\La+1}_3a_2=-a_8,\quad
u^{\La+1}_4a_2=(\La+2)a_6,\allowdisplaybreaks\\ &u^{\La+1}_5
a_2=-a_9-a_{10}+2(\La+2)\pa a_1,\quad u^{\La+1}_6a_2=0,\quad
u^{\La+1}_7a_2=a_{13}-2\pa
a_3,\allowdisplaybreaks\\&u^{\La+1}_8a_2=0,\quad
u^{\La+1}_9a_2=-a_{12}+2\pa a_2,\quad u^{\La+1}_{10}a_2=a_{12}+2(\La+2)\pa a_2,\allowdisplaybreaks\\
&u^{\La+1}_{11}a_2=2(\La+2)\pa a_4-(\La+2)a_{14},\quad
u^{\La+1}_{12}a_2=0,\quad u^{\La+1}_{13}a_2=-2\pa a_8,\allowdisplaybreaks\\
&u^{\La+1}_{14}a_2=2(\La+2)\pa a_6,\quad
u^{\La+1}_{15}a_2=-(\La+2)a_{16}+2(\La+1)\pa a_9-2\pa a_{10},\allowdisplaybreaks\\
&u^{\La+1}_{16}a_2=-2\pa a_{12}.
\end{align*}
Next we find $\C[\pa]$-generators of $N_3^{E_0}$. Similarly
$\{u^{\La+1}_ia_3|1\le i\le 16\}$ generates $N_3^{E_0}$ over
$\C[\pa]$:
\begin{align*}
&u^{\La+1}_1a_3=a_3,\quad u^{\La+1}_2a_3=a_8,\quad
u^{\La+1}_3a_3=0,\quad
u^{\La+1}_4a_3=(\La+1)a_9-a_{10},\allowdisplaybreaks\\
&u^{\La+1}_5 a_3=(\La+2)a_7,\quad u^{\La+1}_6a_3=a_{12},\quad
u^{\La+1}_7a_3=0,\quad
u^{\La+1}_8a_3=0,\allowdisplaybreaks\\
&u^{\La+1}_9a_3=a_{13},\quad u^{\La+1}_{10}a_3=(\La+2)a_{13},\quad
u^{\La+1}_{11}a_3=-(\La+2)a_{15},\allowdisplaybreaks\\
&u^{\La+1}_{12}a_3=-a_{16},\quad u^{\La+1}_{13}a_3=0,\quad
u^{\La+1}_{14}a_3=(\La+2)a_{16},\quad u^{\La+1}_{15}a_3=0,\allowdisplaybreaks\\
&u^{\La+1}_{16}a_3=0.
\end{align*}
It follows that $N^{E_0}$ is generated over $\C[\pa]$ by the
following set
\begin{equation*} S^\La=\{a_2,a_3,a_6,a_7,a_8,a_9-2\pa
a_1,a_{10}-2(\La+1)\pa a_1,a_{12},a_{13},a_{14}-2\pa
a_4,a_{15},a_{16}\}.
\end{equation*}

Suppose that $\La\ge 2$. From the description of $S^\La$ we see
that $\{a_1,a_4,a_5,a_{11}\}$ is a $\C[\pa]$-basis for the
$E_0$-invariants of the quotient
$M_{\NS^4_+}(\frac{\La}{2},\La)/N$. Since $a_1$ has $H_0$-weight
$\La$, it generates a copy of the irreducible $sl_2$-module of
dimension $\La+1$. Now $a_4$ and $a_5$ both have weight $\La-1$,
and so they generate two copies of the irreducible $sl_2$-module
of dimension $\La$. Finally $a_{11}$ has weight $\La-2$, and so it
generates a copy of the irreducible $sl_2$-module of dimension
$\La-1$. Thus $M_{\NS^4_+}(\frac{\La}{2},\La)/N$ is a free
$\C[\pa]$-module of rank $(\La+1)+2\La+(\La-1)=4\La$. So we need
to show that $M_{\NS^4_+}(\frac{\La}{2},\La)/N$ is irreducible.

As in \secref{N=3} $L_n$, $n\ge -1$, together with $E_0,H_0,F_0$
generate a copy of $(\V_+\oplus sl_2)$, which thus allow us to
study the $(\V_+\oplus sl_2)$-module structure of
$M_{\NS^4_+}(\frac{\La}{2},\La)/N$. By parity consideration
$M_{\NS^4_+}(\frac{\La}{2},\La)/N$ is a direct sum of two
modules, namely $(M_{\NS^4_+}(\frac{\La}{2},\La)/N)_{\bar{0}}
=\C[\pa]V_1+\C[\pa]V_{11}$ and
$(M_{\NS^4_+}(\frac{\La}{2},\La)/N)_{\bar{1}}=\C[\pa]V_4+\C[\pa]V_5$,
where $V_i$ is the irreducible $sl_2$-module generated by $a_i$.
We can easily check that $L_n$, for $n\ge 1$, annihilates the
vectors $a_1,a_4,a_5,a_{11}$. (Again the only non-trivial part is
to check that $L_1 a_{11}=0$.) Thus
$M_{\NS^4_+}(\frac{\La}{2},\La)/N$ as a $\V_+\oplus sl_2$-module
is a direct sum of the following four irreducible modules:
$\C[\pa]V_1\cong L_{\V_+}(\frac{\La}{2})\boxtimes U^\La$,
$\C[\pa]V_4\cong L_{\V_+}(\frac{\La+1}{2})\boxtimes U^{\La-1}$,
$\C[\pa]V_5\cong L_{\V_+}(\frac{\La+1}{2})\boxtimes U^{\La-1}$ and
$\C[\pa]V_{11}\cong L_{\V_+}(\frac{\La+2}{2})\boxtimes U^{\La-2}$,
where as usual $U^\mu$ is the irreducible $sl_2$-module of highest
weight $\mu$. Note that, contrary to the $K(1,3)_+$ case, the odd
part here is a sum of two isomorphic modules. To conclude that
$M_{\NS^4_+}(\frac{\La}{2},\La)/N$ is irreducible, we show again
that one may go from each irreducible $\V_+\oplus sl_2$-component
to the irreducible component containing the $\G$-highest weight
vectors. But this follows from the following computation.
\begin{align}
&G^{++}_{\frac{1}{2}}(\alpha a_4+\beta
a_5)=\beta\La(2\La+2)a_1, \quad\alpha,\beta\in\C,\label{N=4aux1}\\
&G^{--}_{\frac{1}{2}}a_4=(2\La+2)F_0 a_1,\quad
E_1a_{11}=\La(\La-1)(2\La+2) a_1.\label{N=4aux2}
\end{align}

Now if $\La=1$ the vector $a_{11}$ is zero. Therefore the quotient
$M_{\NS^4_+}(\frac{\La}{4},\La)/N=\C[\pa]V_1\oplus\C[\pa]V_4\oplus
\C[\pa]V_5$. But then \eqnref{N=4aux1} and the first identity in
\eqnref{N=4aux2} show that $M_{\NS^4_+}(\frac{\La}{4},\La)/N$ is
irreducible.  The rank of $L_{\NS^4_+}(\frac{\La}{4},\La)$ is then
$(\La+1)+2\La$, which equals to $4\La$, in the case $\La=1$.

Finally when $\La=0$, the vectors
$a_4=a_5=a_{10}=a_{11}=a_{14}=a_{15}=0$ so that $S^\La$ reduces to
$\{a_2,a_3,a_6,a_7,a_8,a_9,\pa a_1,a_{12},a_{13},a_{16}\}$. Hence
$M_{\NS^4_+}(0,0)/N=\C a_1$ is the trivial module and so has rank
$0$.
\end{proof}

\begin{prop}
Suppose that $2\Delta+\La+2=0$. Then $L_{\NS^4_+}(\Delta,\La)$ is
a free $\C[\pa]$-module of rank $4\La+8$.
\end{prop}

\begin{proof}
By \propref{singularN=4} $a_4$ and $a_5$ are singular vectors of
$M_{\NS^4_+}(-\frac{\La+2}{2},\La)$ in the case $\La\ge 1$.

Assume first that $\La\ge 3$. Let $N_4$ and $N_5$ be the
$\G$-submodules generated by $a_4$ and $a_5$, respectively.   We
form the $\G$-submodule $N=N_4+N_5$ and consider $N^{E_0}$. The
set $\{u^{\La-1}_ia_4|1\le i\le 16\}$ is a set of
$\C[\pa]$-generators for $N_4^{E_0}$, since $a_4$ has
$H_0$-weight $\La-1$.  We have
\begin{align*}
&u^{\La-1}_1a_4=a_4,\quad u^{\La-1}_2a_4=-\La a_6,\quad
u^{\La-1}_3a_4=2\La\pa a_1-\La a_9+a_{10} ,\quad
u^{\La-1}_4a_4=0,\allowdisplaybreaks\\ &u^{\La-1}_5
a_4=-a_{11},\quad u^{\La-1}_6a_4=0,\quad
u^{\La-1}_7a_4=-a_{15}+2\pa a_5 ,\quad u^{\La-1}_8a_4=2\La\pa
a_2+\La a_{12},\allowdisplaybreaks\\& u^{\La-1}_9a_4=a_{14}+2\pa
a_4,\quad u^{\La-1}_{10}a_4=(\La-1)a_{14},\quad
u^{\La-1}_{11}a_4=0,\quad u^{\La-1}_{12}a_4=2\La\pa
a_6,\allowdisplaybreaks\\
&u^{\La-1}_{13}a_4=-\La a_{16}+2\pa a_{10} ,\quad
u^{\La-1}_{14}a_4=0,\quad u^{\La-1}_{15}a_4=-2\pa a_{11},\quad
u^{\La-1}_{16}a_4=2\pa a_{14}.
\end{align*}
Similarly, the following is a set of $\C[\pa]$-generators for
$N_5^{E_0}$.
\begin{align*}
&u^{\La-1}_1a_5=a_5,\quad u^{\La-1}_2a_5=a_{10},\quad
u^{\La-1}_3a_5=-a_7,\quad u^{\La-1}_4a_5=a_{11},\quad u^{\La-1}_5
a_5=0,\allowdisplaybreaks\\
&u^{\La-1}_6a_5=a_{14},\quad u^{\La-1}_7a_5=0,\quad
u^{\La-1}_8a_5=-\La a_{13},\quad u^{\La-1}_9a_5=a_{15},\quad
u^{\La-1}_{10}a_5=0,
\allowdisplaybreaks\\
&u^{\La-1}_{11}a_5=0,\quad u^{\La-1}_{12}a_5=-a_{16},\quad
u^{\La-1}_{13}a_5=0,\quad u^{\La-1}_{14}a_5=0,\quad
u^{\La-1}_{15}a_5=0,\\
&u^{\La-1}_{16}a_5=0.
\end{align*}
Therefore $$S^\La=\{a_4,a_5,a_6,a_7,a_9-2\pa
a_1,a_{10},a_{11},a_{12}+2\pa a_2,a_{13},a_{14},a_{15},a_{16}\}$$
is a set of $\C[\pa]$-generators for $N^{E_0}$, which implies that
$\{a_1,a_2,a_3,a_8\}$ is a $\C[\pa]$-basis for
$(M_{\NS^4_+}(-\frac{\La+2}{2},\La)/N)^{E_0}$ in the case when
$\La\ge 3$.

In the case when $\La=2$ the set
$\{u^{\La-1}_ia_4,u^{\La-1}_ia_5|1\le i\le 16,i\not=11\}$ generate
$N^{E_0}$. But $u^{\La-1}_{11}a_4=u^{\La-1}_{11}a_5=0$, and hence
$\{a_1,a_2,a_3,a_8\}$ is also a $\C[\pa]$-basis for
$(M_{\NS^4_+}(-\frac{\La+2}{2},\La)/N)^{E_0}$ in this case as
well.

In the case when $\La=1$ we note that $a_{11}=0$ and
$\{u^{\La-1}_ia_4,u^{\La-1}_ia_5|1\le i\le
16,i\not=4,5,10,11,14,15\}$ generates $N^{E_0}$ over $\C[\pa]$.
From the formulas above one sees that a set of
$\C[\pa]$-generators  for $N^{E_0}$ is given by the set $S^\La$
above, but with $a_{11}$ removed. Hence the quotient module is
again generated freely over $\C[\pa]$ by $\{a_1,a_2,a_3,a_8\}$.

Hence in the case when $\La\ge 1$ the quotient module
$(M_{\NS^4_+}(-\frac{\La+2}{2},\La)/N)^{E_0}$ is generated freely
over $\C[\pa]$ by $\{a_1,a_2,a_3,a_8\}$.  Now $a_1$ has
$H_0$-weight $\La$, $a_2$ and $a_3$ both have $H_0$-weight
$\La+1$, and $a_8$ has $H_0$-weight $\La+2$. Therefore
$M_{\NS^4_+}(-\frac{\La+2}{2},\La)/N$ has rank
$(\La+1)+2(\La+2)+(\La+3)=4\La+8$ over $\C[\pa]$. So it remains
to show that $M_{\NS^4_+}(-\frac{\La+2}{2},\La)/N$ is irreducible.

We again study $M_{\NS^4_+}(-\frac{\La+2}{2},\La)/N$ as a
$\V_+\oplus sl_2$-module.  It is easy to check that $L_n$, $n\ge
1$, annihilates $a_1,a_2,a_3,a_8$ and hence
$M_{\NS^4_+}(-\frac{\La+2}{2},\La)/N$ is a direct sum of the
following four irreducible $\V_+\oplus sl_2$-modules:
$\C[\pa]V_1\cong L_{\V_+}(-\frac{\La+2}{2})\boxtimes U^\La$,
$\C[\pa]V_2\cong L_{\V_+}(-\frac{\La+1}{2})\boxtimes U^{\La+1}$,
$\C[\pa]V_3\cong L_{\V_+}(-\frac{\La+1}{2})\boxtimes U^{\La+1}$
and $\C[\pa]V_8\cong L_{\V_+}(-\frac{\La}{2})\boxtimes U^{\La+2}$,
where $V_i$ is the $sl_2$-submodule generated by the vector
$a_i$. Again $\C[\pa]V_2\cong\C[\pa]V_3$ as $\V_+\oplus
sl_2$-modules. Now we compute
\begin{align*}
&G^{--}_{\half}(\alpha a_2+\beta
a_3)=2\alpha(\La+1)a_1,\quad\forall \alpha,\beta\in\C,\\
&G^{-+}_\half a_3=-2(\La+1)a_1,\quad F_1 a_8=-2(\La+1)a_1.
\end{align*}
Therefore $M_{\NS^4_+}(-\frac{\La+2}{2},\La)/N$ is irreducible.

Now consider the case of $\La=0$.  By \propref{singularN=4} $a_6$,
$a_7$ and $a_9-2\pa a_1$ are singular vectors inside
$M_{\NS^4_+}(-1,0)$. Let $N_6$, $N_7$ and $N_9$ be the
$\G$-submodules generated by $a_6$, $a_7$ and $a_9-2\pa a_1$,
respectively, and put $N=N_6+N_7+N_9$. We note that $a_6$, $a_7$
and $a_9-2\pa a_1$ have $H_0$-weight $0$, hence $N_6^{E_0}$ is
generated over $\C[\pa]$ by $\{u^{\La}_ia_6|1\le i\le
16,i\not=4,5,10,11,14,15\}$ and similarly for $N_7^{E_0}$ and
$N_9^{E_0}$.  We first compute a set of $\C[\pa]$-generators for
$N_6^{E_0}$.
\begin{align*}
&u^{\La}_1a_6=a_6,\quad u^{\La}_2a_6=0,\quad u^{\La}_3a_6=2\pa a_2+a_{12},\quad u^{\La}_6a_6=0,\\
&u^{\La}_7a_6=4\pa^2 a_1-2\pa a_9+a_{16},\quad u^{\La}_8a_6=0,\quad u^{\La}_9a_6=4\pa a_6,\\
&u^{\La}_{12}a_6=0,\quad u^{\La}_{13}a_6=4\pa^2 a_2+2\pa
a_{12},\quad u^{\La}_{16}a_6=4\pa^2 a_6.
\end{align*}
A set of $\C[\pa]$-generators for $N_7^{E_0}$ is given as follows.
\begin{align*}
&u^{\La}_1a_7=a_7,\quad u^\La_2a_7=a_{13},\quad
u^{\La}_3a_7=0,\quad u^{\La}_6a_7=a_{16},\quad u^{\La}_7a_7=0, \\
&u^{\La}_8a_7=0,\quad u^{\La}_9a_7=0,\quad u^{\La}_{12}a_7=0,\quad
u^{\La}_{13}a_7=0,\quad u^{\La}_{16}a_7=0.
\end{align*}
Finally we have the following set of $\C[\pa]$-generators for
$N_9^{E_0}$.
\begin{align*}
&u^{\La}_1(a_9-2\pa a_1)=a_9-2\pa a_1,\quad u^{\La}_2(a_9-2\pa
a_1)=-a_{12}-2\pa a_2,\\
&u^{\La}_3(a_9-2\pa a_1)=a_{13},\quad u^{\La}_6(a_9-2\pa
a_1)=-2\pa a_6,\quad u^{\La}_7(a_9-2\pa a_1)=2\pa a_7,\\
&u^{\La}_8(a_9-2\pa
a_1)=0,\quad u^{\La}_9(a_9-2\pa a_1)=2a_{16},\quad u^{\La}_{12}(a_9-2\pa a_1)=0,\\
&u^{\La}_{13}(a_9-2\pa a_1)=2\pa a_{13},\quad
u^{\La}_{16}(a_9-2\pa a_1)=2\pa a_{16}.
\end{align*}
From this it follows that $\{a_6,a_7,a_9-2\pa a_1,a_{12}+2\pa
a_2,a_{13},a_{16}\}$ generate $N^{E_0}$ over $\C[\pa]$. But
$a_4=a_5=a_{10}=a_{11}=a_{14}=a_{15}=0$, and thus
$(M_{\NS^4_+}(-1,0)/N)^{E_0}$ is generated over $\C[\pa]$ by the
vectors $a_1$, $a_2$, $a_3$ and $a_8$, which takes us back to the
case when $\La\ge 1$, except that here $\C[\pa]V_8$ is not
irreducible.  It contains a unique irreducible submodule
isomorphic to $L_{\V_+}(1)\boxtimes U^2$ generated by $\pa a_8$.
But then the above calculation plus the fact that $$F_2 \pa
a_8=-4(\La+1)a_1$$ show that $M_{\NS^4_+}(-1,0)/N$ is irreducible
of rank $8$.
\end{proof}

\begin{thm}
The modules $L_{\NS^4_+}(\Delta,\La)$, for $\Delta\in\C$ and
$\La\in\Z_+$, form a complete list of non-isomorphic finite (over
$\C[\pa]$) irreducible $SK(1,4)_+$-modules.  Furthermore
$L_{\NS^4_+}(\Delta,\La)$ as a $\C[\pa]$-module has rank
\begin{itemize}
\item[i.] $4\La$, in the case $2\Delta-\La=0$,
\item[ii.] $4\La+8$, in the case $2\Delta+\La+2=0$.
\item[iii.] $16\La+16$, in all other cases.
\end{itemize}
Furthermore the $\C[\pa]$-rank of
$L_{\NS^4_+}(\Delta,\La)_{\bar{0}}$ equals the $\C[\pa]$-rank of
$L_{\NS^4_+}(\Delta,\La)_{\bar{1}}$ in all cases.
\end{thm}

\begin{rem}
Translating the above theorem into the languages of modules over
conformal algebras and of conformal modules is again
straightforward.  We therefore obtain that all finite irreducible
modules over the ``small'' $N=4$ conformal superalgebra are of the
form $L_{\NS^4}(\alpha,\Delta,\La)$, where $\alpha,\Delta\in\C$
and $\La\in\Z_+$.  The definition of these modules and also the
action of the conformal superalgebra on them are easily gotten
from our explicit description of a $\C[\pa]$-basis in this
section and hence omitted, as to reproduce them would take up
quite a significant portion of space. Again we only note that the
adjoint module is isomorphic to $L_{\NS^4}(0,1,2)$.
\end{rem}

\section{Finite irreducible Modules over the ``big'' $N=4$ conformal superalgebra}\label{BN=4}

In this section we give a classification of finite irreducible
conformal modules over the contact superalgebra $K(1,4)$, also
known as the ``big'' $N=4$ superconformal algebra. Our approach
is based on our results obtained in \secref{N=4}.

Recall from \secref{N=4} that $L_n^\beta$, $X^\beta_n$ and
$x^\beta_r$, where $X=H,E,F$, $x=G^{++}$,
$G^{-+}$,$G^{+-}$,$G^{--}$, $n\in\Z$, $r\in\frac{1}{2}+\Z$ and
the fixed number $\beta$ is either $1$ or $-1$, provide a basis
for a copy of $SK(1,4)$ inside $K(1,4)$. In this section it will
be convenient to distinguish these two copies.  We therefore
denote the copy obtained by setting $\beta=1$ simply by
$SK(1,4)$, while the copy obtained by setting $\beta=-1$ by
$\overline{SK}(1,4)$. It is easy to see from our formulas that
$K(1,4)=SK(1,4)+\overline{SK}(1,4)$. Similarly we distinguish the
basis elements of $SK(1,4)$ and $\overline{SK}(1,4)$ as follows.
The generators inside $SK(1,4)$ will be denoted by $L_n,X_n,x_r$,
while generators inside $\overline{SK}(1,4)$ will be denoted by
$\ov{L}_n,\ov{X}_n,\ov{x}_r$, where again $X=H,E,F$,
$x=G^{++}$,$G^{-+}$,$G^{+-}$,$G^{--}$. Of course we have
$x_{-\half}=\ov{x}_{-\half}$, $L_{-1}=\ov{L}_{-1}$ and
$L_0=\ov{L}_0$.

\begin{rem}\label{symmetry}
The map $\phi:SK(1,4)\rightarrow \ov{SK}(1,4)$ defined by
$\phi(L_{n})=\ov{L}_n$, $\phi(X_{n})=\ov{X}_n$,
$\phi(G^{++}_{r})=\ov{G}^{++}_r$,
$\phi(G^{+-}_{r})=\ov{G}^{-+}_r$,
$\phi(G^{-+}_{r})=\ov{G}^{+-}_r$ and
$\phi(G^{--}_{r})=\ov{G}^{--}_r$, where $n\in\Z$ and
$r\in\half+\Z$ is an isomorphism of Lie superalgebras.  Thus all
formulas in \secref{N=4} with $\phi(L_n)$, $\phi(X_n)$ and
$\phi(x_r)$ replacing $L_n$, $X_n$ and $x_r$, respectively,
remain valid.
\end{rem}

Let $\G=K(1,4)_+$ be the annihilation subalgebra of $K(1,4)$ so
that we have $\G=SK(1,4)_++\overline{SK}(1,4)_+$, the sum of the
corresponding annihilation subalgebras. We have as before
$\G=\bigoplus_{j\ge -1 }\G_{j}$, where $j\in\half+\Z$. Furthermore
$\G_-=SK(1,4)_-=\overline{SK}(1,4)_-$ and $\G_0=\C L_0\oplus
sl_2\oplus\overline{sl}_2\cong cso_4$, where $sl_2$ and
$\overline{sl}_2$ denote two copies of the Lie algebra $sl_2$,
generated by $H_0,E_0,F_0$ and $\ov{H}_0,\ov{E}_0,\ov{F}_0$,
respectively.

Let $U^{\Delta,\La,\ov{\La}}$ be the finite-dimensional
irreducible $sl_2\oplus\ov{sl}_2$-module of highest weight
$(\La,\ov{\La})\in\Z_+\times\Z_+$ on which $L_0$ acts as the
scalar $\Delta\in\C$ and let $v_{\Delta,\La,\ov{\La}}$ denote a
highest weight vector in $U^{\Delta,\La,\ov{\La}}$ so that
$H_0v_{\Delta,\La,\ov{\La}}=\La v_{\Delta,\La,\ov{\La}}$,
$\ov{H}_0v_{\Delta,\La,\ov{\La}}=\ov{\La} v_{\Delta,\La,\ov{\La}}$
and $L_0v_{\Delta,\La,\ov{\La}}=\Delta v_{\Delta,\La,\ov{\La}}$.
Regarding $U^{\Delta,\La,\ov{\La}}$ as a module over
$\Li_0=\bigoplus_{j\ge 0}\G_j$ it follows from \thmref{key} that
every finite irreducible $\G$-module is a quotient of
$M_{\SN^4_+}(\Delta,\La,\ov{\La})={\rm Ind}_{\Li_0}^\G
U^{\Delta,\La,\ov{\La}}$.  The unique irreducible quotient will
be denoted by $L_{\SN^4_+}(\Delta,\La,\ov{\La})$.

Now $M_{\SN^4_+}(\Delta,\La,\ov{\La})$ is a completely reducible
$\G_0$-module, and the subspace of $\C
E_0\oplus\C\ov{E}_0$-invariants, denoted by
$M_{\SN^4_+}(\Delta,\La,\ov{\La})^{E_0,\ov{E}_0}$, is a free
$\C[\pa]$-submodule of $M_{\SN^4_+}(\Delta,\La,\ov{\La})$. We
write down explicit formulas for a $\C[\pa]$-basis for the space
$M_{\SN^4_+}(\Delta,\La,\ov{\La})^{E_0,\ov{E}_0}$, which in the
case when $\La,\ov{\La}\ge 2$ is as follows:
\begin{align*}
&b_1=v_{\Delta,\La,\ov{\La}},\quad
b_2=G^{++}_{-\frac{1}{2}}v_{\Delta,\La,\ov{\La}},\allowdisplaybreaks\\
&b_3=(\La G^{-+}_{-\half}-G^{++}_{-\half}
F_0)v_{\Delta,\La,\ov{\La}},\quad b_4=(\ov{\La}
G^{+-}_{-\half}-G^{++}_{-\half}
\ov{F}_0)v_{\Delta,\La,\ov{\La}},\allowdisplaybreaks\\
&b_5=(\La\ov{\La}G^{--}_{-\half}-\ov{\La}G^{+-}_{-\half}F_0-\La
G^{-+}_{-\half}\ov{F}_0 +G^{++}_{-\half}F_0\ov{F}_0)
v_{\Delta,\La,\ov{\La}},\allowdisplaybreaks\\
&b_6=G^{++}_{-\frac{1}{2}}G^{+-}_{-\frac{1}{2}}v_{\Delta,\La,\ov{\La}},\quad
b_7=\big(\La(G^{-+}_{-\frac{1}{2}}G^{+-}_{-\frac{1}{2}}+G^{++}_{-\half}G^{--}_{-\half})-
2G^{++}_{-\half}G^{+-}_{-\half}F_0\big)v_{\Delta,\La,\ov{\La}},\allowdisplaybreaks\\
&b_8=((\La-1)(-\La G^{-+}_{-\frac{1}{2}}G^{--}_{-\frac{1}{2}}+
G^{-+}_{-\frac{1}{2}}G^{+-}_{-\frac{1}{2}}
F_0+G^{++}_{-\frac{1}{2}}G^{--}_{-\frac{1}{2}} F_0)
-G^{++}_{-\frac{1}{2}}
G^{+-}_{-\frac{1}{2}}F_0^2)v_{\Delta,\La,\ov{\La}},
\allowdisplaybreaks\\
&b_9=G^{++}_{-\half}G^{-+}_{-\half}v_{\Delta,\La,\ov{\La}},\quad
b_{10}=\big(\ov{\La} (G^{++}_{-\frac{1}{2}}G^{--}_{-\frac{1}{2}}-
G^{-+}_{-\frac{1}{2}}G^{+-}_{-\frac{1}{2}})-
2G^{++}_{-\half}G^{-+}_{-\half}\ov{F_0}\big)v_{\Delta,\La,\ov{\La}},\allowdisplaybreaks\\
&b_{11}=\big((\ov{\La}-1)(-\ov{\La}
G^{+-}_{-\frac{1}{2}}G^{--}_{-\frac{1}{2}}+
G^{+-}_{-\frac{1}{2}}G^{-+}_{-\frac{1}{2}}
\ov{F}_0+G^{++}_{-\frac{1}{2}}G^{--}_{-\frac{1}{2}} \ov{F}_0)
-G^{++}_{-\frac{1}{2}} G^{-+}_{-\frac{1}{2}}\ov{F}_0^2\big) v_{\Delta,\La,\ov{\La}},\allowdisplaybreaks\\
&b_{12}=G^{++}_{-\frac{1}{2}}G^{+-}_{-\frac{1}{2}}G^{-+}_{-\frac{1}{2}}
v_{\Delta,\La,\ov{\La}},\quad b_{13}=(\La
G^{++}_{-\half}G^{-+}_{-\half}G^{--}_{-\half}
-G^{++}_{-\half}G^{-+}_{-\half}G^{+-}_{-\half})v_{\Delta,\La,\ov{\La}},
\allowdisplaybreaks\\
&b_{14}=(\ov{\La} G^{++}_{-\half}G^{+-}_{-\half}G^{--}_{-\half}
-G^{++}_{-\half}G^{+-}_{-\half}G^{-+}_{-\half})v_{\Delta,\La,\ov{\La}},
\allowdisplaybreaks\\
&b_{15}=\big(\La
G^{-+}_{-\half}G^{--}_{-\half}(\ov{\La}G^{+-}_{-\half}-
G^{++}_{-\half}\ov{F}_0)+
\ov{\La}G^{++}_{-\half}G^{+-}_{-\half}(G^{--}_{-\half}F_0-
G^{-+}F_0\ov{F}_0)\big) v_{\Delta,\La,\ov{\La}}
,\allowdisplaybreaks\\
&b_{16}=\big(G^{++}_{-\frac{1}{2}}G^{+-}_{-\frac{1}{2}}
G^{-+}_{-\frac{1}{2}}G^{--}_{-\frac{1}{2}}-
\pa(G^{-+}_{-\half}G^{+-}_{-\half}
+G^{++}_{-\half}G^{--}_{-\half})\big)v_{\Delta,\La,\ov{\La}}.
\end{align*}
In the case when $\La=\ov{\La}=1$ (respectively $\La=\ov{\La}=0$)
we have $b_8=b_{11}=0$ (respectively
$b_3=b_4=b_5=b_7=b_8=b_{10}=b_{11}=b_{13}=b_{14}=b_{15}=0$), thus
giving us $14$ (respectively $6$) generators.  Other cases are
easily described as well, however, we will not need them because
of \propref{BN=4generic1} below.  Thus we will omit them.

We will, as before, denote the coefficient of
$v_{\Delta,\La,\ov{\La}}$ in $b_i$ by $u_i^{\La,\ov{\La}}$ for
$1\le i\le 16$. In the case when $\La=\ov{\La}$, which is the
only case we will be concerned with in what follows, we simply
write $u_i^\La$ for $u_i^{\La,\La}$ and also $v_{\Delta,\La}$ for
$v_{\Delta,\La,\La}$.

\begin{prop}\label{BN=4generic1}
If $M_{\SN^4_+}(\Delta,\La,\ov{\La})$ is a reducible $\G$-module,
then either $2\Delta-\La=2\Delta-\ov{\La}=0$ or else
$2\Delta+\La+2=2\Delta+\ov{\La}+2=0$.  In particular if
$\La\not=\ov{\La}$, then $M_{\SN^4_+}(\Delta,\La,\ov{\La})$ is
irreducible.
\end{prop}

\begin{proof}
As a module over $SK(1,4)_+$ we have
$M_{\SN^4_+}(\Delta,\La,\ov{\La})=U(\G_-)\otimes
U^{\Delta,\La,\ov{\La}}$ is a direct sum of $\ov{\La}+1$ copies of
$M_{\NS^4_+}(\Delta,\La)$, generated by the highest weight vectors
$\ov{F}_0^{j}v_{\Delta,\La,\ov{\La}}$, where $0\le j\le
\ov{\La}$.  Since the $\ov{H}_0$-weights of the
$\ov{F}_0^{j}v_{\Delta,\La,\ov{\La}}$'s are all distinct for
distinct $j$'s it follows that these modules as
$SK(1,4)_+\rtimes\C\ov{H}_0$-modules are all non-isomorphic.
Therefore if $M_{\NS^4_+}(\Delta,\La)$ is irreducible over
$SK(1,4)_+$, then $M_{\SN^4_+}(\Delta,\La,\ov{\La})$ is
irreducible over $\G$. From this and \corref{N=4generic} we thus
conclude that in the case when $\Delta-2\La\not=0$ and
$\Delta+2\La+2\not=0$ the $\G$-module
$M_{\SN^4_+}(\Delta,\La,\ov{\La})$ is irreducible.

By symmetry we conclude that if $\Delta-2\ov{\La}\not=0$ and
$\Delta+2\ov{\La}+2\not=0$, then
$M_{\SN^4_+}(\Delta,\La,\ov{\La})$ is irreducible over $\G$ as
well.

Therefore $M_{\SN^4_+}(\Delta,\La,\ov{\La})$ is possibly
reducible only if both $\La$ and $\ov{\La}$ satisfy one of the
two linear equations $\Delta-2x=0$ and $\Delta+2x+2=0$.  But the
case $\Delta-2\La=0$ and $\Delta+2\ov{\La}+2=0$ is not possible,
since both $\La$ and $\ov{\La}$ are non-negative integers. By the
same token $\Delta-2\ov{\La}=0$ and $\Delta+2\La+2=0$ is not
possible, either. Hence either we have $\Delta-2\La=0$ and
$\Delta-2\ov{\La}=0$ or else $\Delta+2\La+2=0$ and
$\Delta+2\ov{\La}+2=0$.  In either case we must have
$\La=\ov{\La}$.
\end{proof}

The next step is to analyze proper singular vectors inside
$M_{\SN^4}(\Delta,\La,\ov{\La})$. (The definitions of singular
vectors and proper singular vectors of $\G$ are of course
analogous.) By \propref{BN=4generic1} proper singular vectors
exist only if $\La=\ov{\La}$ with either $2\Delta+\La=0$ or
$2\Delta+\La+2=0$.

\begin{prop}\label{singularBN=4}
A complete list of proper singular vectors inside
$M_{\SN^4_+}(\Delta,\La,\La)$ is given by:
\begin{itemize}
\item[i.] $\alpha b_2$, $\alpha\not=0$, in the case $2\Delta-\La=0$.
\end{itemize}
\begin{itemize}
\item[ii.] $\alpha b_5$, $\alpha\not=0$, in the case
$2\Delta+\La+2=0$ and $\La\ge 1$.
\end{itemize}
\end{prop}

\begin{proof}
Since as a $SK(1,4)_+$-module $M_{\SN^4}(\Delta,\La,\La)$ is a
direct sum of $\La+1$ copies of $M_{\NS^4}(\Delta,\La)$ we obtain
a description of the vector space spanned by all proper
$SK(1,4)_+$-singular vectors by virtue of \propref{singularN=4}.
But as a $\ov{SK}(1,4)_+$-module $M_{\SN^4}(\Delta,\La,\La)$ is
also a direct sum of $\La+1$ copies of $M_{\NS^4}(\Delta,\La)$,
from which we obtain similarly a description of the vector space
spanned by all proper $\ov{SK}(1,4)_+$-singular vectors (see
\remref{symmetry}).  The intersection of these two spaces is the
space of proper singular vectors.

In the case when $2\Delta-\La=0$ it follows from
\propref{singularN=4} that the space of proper
$SK(1,4)_+$-singular vectors is spanned by
$G^{++}_{-\half}\ov{F}_0^jv_{\Delta,\La}$,
$G^{+-}_{-\half}\ov{F}_0^jv_{\Delta,\La}$ and
$G^{++}_{-\half}G^{+-}_{-\half}\ov{F}_0^jv_{\Delta,\La}$, for
$0\le j\le \La$.  On the other hand the space of proper
$\ov{SK}(1,4)_+$-singular vectors is spanned by
$G^{++}_{-\half}F_0^jv_{\Delta,\La}$,
$G^{-+}_{-\half}F_0^jv_{\Delta,\La}$ and
$G^{++}_{-\half}G^{-+}_{-\half}F_0^jv_{\Delta,\La}$, for $0\le
j\le \La$. It is not hard to see that the intersection of these
two spaces is the one-dimensional space spanned by
$G^{++}_{-\half}v_{\Delta,\La}$, which is $b_2$.

Other cases are analogous and so we omit the details.
\end{proof}

\begin{prop}
Suppose that $2\Delta-\La=0$. Then $L_{\SN^4_+}(\Delta,\La,\La)$
is a free $\C[\pa]$-module of rank $8\La(\La+1)$.
\end{prop}

\begin{proof}
By \propref{singularBN=4} $b_2$ is a singular vector in
$M_{\SN^4_+}(\frac{\La}{2},\La,\La)$. Consider $N$, the
$\G$-submodule generated by $b_2$.  Then we have $N=U(\G_{-})V$,
where $V$ is the irreducible $sl_2\oplus\ov{sl}_2$-submodule
generated by $b_2$. Let us compute the space $N^{E_0,\ov{E}_0}$,
the space of $(\C E_0\oplus\C\ov{E}_0)$-invariants inside $N$.
Since the $(H_0,\ov{H}_0)$-weight of $b_2$ is $(\La+1,\La+1)$, we
know that $N^{E_0,\ov{E}_0}$ is a free $\C[\pa]$-module generated
over $\C[\pa]$ by $\{u^{\La+1}_ib_2|1\le i\le 16\}$. We have
\begin{align*}
&u^{\La+1}_1b_2=b_2,\quad u^{\La+1}_2b_2=0,\quad
u^{\La+1}_3b_2=-(\La+2)b_9,\quad
u^{\La+1}_4b_2=-(\La+2)b_6,\allowdisplaybreaks\\ &u^{\La+1}_5
b_2=-{\frac{\La+2}{2}}(b_7+b_{10}+4(\La+1)\pa b_1),\quad
u^{\La+1}_6b_2=0,\quad
u^{\La+1}_7b_2=-(\La+3)b_{12},\allowdisplaybreaks\\&u^{\La+1}_8b_2=-(\La+2)(b_{13}-2\pa
b_3) ,\quad
u^{\La+1}_9b_2=0,\quad u^{\La+1}_{10}b_2=(\La+3)b_{12}-4(\La+1)\pa b_2,\allowdisplaybreaks\\
&u^{\La+1}_{11}b_2=-(\La+2)(b_{14}-2\pa b_4),\quad
u^{\La+1}_{12}b_2=0,\quad u^{\La+1}_{13}b_2=-2(\La+2)\pa b_9,\allowdisplaybreaks\\
&u^{\La+1}_{14}b_2=-2(\La+2)\pa b_6,\quad
u^{\La+1}_{15}b_2=-4(\La+1)\pa^2 b_1-(\La+2)^2b_{16}+(\La+2)\pa b_7,\allowdisplaybreaks\\
&u^{\La+1}_{16}b_2=-4\pa b_{12}.
\end{align*}
It follows that $N^{E_0,\ov{E}_0}$ is generated over $\C[\pa]$ by
the set
\begin{align*}
S^\La=\{&b_2,b_6,b_7+b_{10}+4(\La+2)\pa b_1,b_9,
b_{12},b_{13}-2\pa b_3,b_{14}-2\pa
b_4,\\&b_{16}-(\frac{1}{\La+2})\pa b_7-
2\frac{(\La+1)}{(\La+2)^2}\pa^2 b_1\}.
\end{align*}
In the case when $\La\ge 2$ it follows from the description of
$S^\La$ that $\{b_1,b_3,b_4,b_5,b_8,\break b_{10}+2\La\pa b_1
,b_{11},b_{15}\}$ is a $\C[\pa]$-basis for the $(\C
E_0\oplus\C\ov{E}_0)$-invariants of the quotient space
$M_{\SN^4_+}(\frac{\La}{2},\La,\La)/N$. (The choice of
$b_{10}+2\La\pa b_1$ instead of just $b_{10}$ will be explained
later.)

The $(L_0,H_0,\ov{H}_0)$-weights of $b_1$, $b_3$, $b_4$, $b_5$,
$b_8$, $b_{10}+2\La\pa b_1$, $b_{11}$, $b_{15}$ are
$(\Delta,\La,\La)$, $(\Delta+\half,\La-1,\La+1)$,
$(\Delta+\half,\La+1,\La-1)$, $(\Delta+\half,\La-1,\La-1)$,
$(\Delta+1,\La-2,\La)$, $(\Delta+1,\La,\La)$,
$(\Delta+1,\La,\La-2)$, $(\Delta+\frac{3}{2},\La-1,\La-1)$,
respectively.  Hence $M_{\SN^4_+}(\frac{\La}{2},\La,\La)/N$ is a
free $\C[\pa]$-module of rank $8\La(\La+1)$. So we need to show
that $M_{\SN^4_+}(\frac{\La}{2},\La,\La)/N$ is irreducible.

Now $L_n$, $n\ge -1$, together with $E_0,H_0,F_0$ and
$\ov{E}_0,\ov{H}_0,\ov{F}_0$ generate a copy of $(\V_+\oplus
sl_2\oplus\ov{sl}_2)$, which thus allow us to study the
$(\V_+\oplus sl_2\oplus\ov{sl}_2)$-module structure of
$M_{\SN^4_+}(\frac{\La}{2},\La,\La)/N$.  We can easily check that
$L_n$, for $n\ge 1$, annihilates the vectors $b_1$, $b_3$, $b_4$,
$b_5$, $b_8$, $b_{10}+2\La\pa b_1$, $b_{11}$, $b_{15}$. (We want
to point out that $b_{10}$ is not annihilated by $L_n$, for $n\ge
1$, hence the choice of $b_{10}+2\La\pa b_1$.) Thus
$M_{\SN^4_+}(\frac{\La}{2},\La,\La)/N$ as a $(\V_+\oplus
sl_2\oplus\ov{sl}_2)$-module is a direct sum of the following
eight irreducible modules: $\C[\pa]V_1\cong
L_{\V_+}(\frac{\La}{2})\boxtimes U^{\La,\La}$, $\C[\pa]V_3\cong
L_{\V_+}(\frac{\La+1}{2})\boxtimes U^{\La-1,\La+1}$,
$\C[\pa]V_4\cong L_{\V_+}(\frac{\La+1}{2})\boxtimes
U^{\La+1,\La-1}$, $\C[\pa]V_{5}\cong
L_{\V_+}(\frac{\La+1}{2})\boxtimes U^{\La-1,\La-1}$,
$\C[\pa]V_8\cong L_{\V_+}(\frac{\La+2}{2})\boxtimes
U^{\La-2,\La}$, $\C[\pa]V_{10}\cong
L_{\V_+}(\frac{\La+2}{2})\boxtimes U^{\La,\La}$,
$\C[\pa]V_{11}\cong L_{\V_+}(\frac{\La+2}{2})\boxtimes
U^{\La,\La-2}$, $\C[\pa]V_{15}\cong
L_{\V_+}(\frac{\La+3}{2})\boxtimes U^{\La-1,\La-1}$, where $V_i$
is the irreducible $sl_2\oplus\ov{sl}_2$-module generated by
$b_i$, for $i\not=10$, and $V_{10}$ is generated by
$b_{10}+2\La\pa b_1 $, and finally $U^{\mu,\mu'}$ denotes the
irreducible $sl_2\oplus\ov{sl}_2$-module of highest weight
$(\mu,\mu')$. Note that as $(\V_+\oplus
sl_2\oplus\ov{sl}_2)$-modules  they are all non-isomorphic and
thus to show that $M_{\SN^4_+}(\frac{\La}{2},\La,\La)/N$ is
irreducible, it suffices to show that one may send a $(\V_+\oplus
sl_2\oplus\ov{sl}_2)$-highest weight vector in any irreducible
$(\V_+\oplus sl_2\oplus\ov{sl}_2)$-component to the irreducible
component containing the $\G$-highest weight vectors. This
follows from the following computation.
\begin{align*}
&G^{--}_{\half}b_3=2(\La+1)F_0b_1,\quad
\ov{G}^{--}_{\half}b_4=2(\La+1)\ov{F}_0b_1,\\
&G^{++}_{\half}b_5=-2\La^2(\La+1)b_1,\quad
E_1b_{8}=2\La(\La-1)(\La+1) b_1,\\
&\ov{F}_1(b_{10}+2\La\pa b_1)=-2(\La+2)\ov{F}_0b_1,\quad
\ov{E}_1b_{11}=2\La(\La-1)(\La+1) b_1,\\
&\ov{G}^{++}_{\frac{3}{2}}b_{15}=-2\La^2(\La+1)b_1.
\end{align*}

Now if $\La=1$ the vectors $b_8=b_{11}=0$. Therefore
$M_{\SN^4_+}(\frac{\La}{2},\La,\La)/N$ is
$\C[\pa]V_1\oplus\C[\pa]V_3\oplus
\C[\pa]V_4\oplus\C[\pa]V_5\oplus\C[\pa]V_{10}\oplus\C[\pa]V_{15}$.
But then the above calculation also shows that
$M_{\SN^4_+}(\frac{\La}{2},\La,\La)/N$ is irreducible.  The rank
of $L_{\SN^4_+}(\frac{\La}{2},\La,\La)$ is then $4+3+3+1+4+1=16$,
which is equal to $8\La(\La+1)$ in the case when $\La=1$.

Finally when $\La=0$, the vectors
$b_3=b_4=b_5=b_7=b_8=b_{10}=b_{11}=b_{13}=b_{14}=b_{15}=0$ and
$S^\La$ reduces to $\{b_2,b_6,\pa b_1,b_9,b_{12},b_{16}\}$. Hence
$M_{\SN^4_+}(0,0,0)/N=\C b_1$ is the trivial module and so has
rank $0$.
\end{proof}

\begin{prop}
Suppose that $2\Delta+\La+2=0$ and $\La\ge 1$. Then
$L_{\SN^4_+}(\Delta,\La,\La)$ is a free $\C[\pa]$-module of rank
$8(\La+1)(\La+2)$.
\end{prop}

\begin{proof}
By \propref{singularBN=4} $b_5$ is a singular vector in
$M_{\SN^4_+}(-\frac{\La+2}{2},\La,\La)$. Let $N$ be the
$\G$-submodule generated by $b_5$ so that $N=U(\G_{-})V$, where
$V$ is the irreducible $sl_2\oplus\ov{sl}_2$-submodule generated
by $b_5$. Consider $N^{E_0,\ov{E}_0}$, the subspace in $N$ of $\C
E_0\oplus\C\ov{E}_0$-invariants. Now the $(H_0,\ov{H}_0)$-weight
of $b_5$ is $(\La-1,\La-1)$ and so $N^{E_0,\ov{E}_0}$ is a free
$\C[\pa]$-module generated over $\C[\pa]$ by
$\{u^{\La-1}_ib_5|1\le i\le 16\}$. We have
\begin{align*}
&u^{\La-1}_1b_5=b_5,\quad
u^{\La-1}_2b_5=\frac{1}{2}(b_7+b_{10}),\quad u^{\La-1}_3b_5=-\La
b_8,\quad u^{\La-1}_4b_5=-\La b_{11},\allowdisplaybreaks\\
&u^{\La-1}_5 b_5=0,\quad u^{\La-1}_6b_5=\La b_{14},\quad
u^{\La-1}_7b_5=-(\La-1)b_{15},\quad
u^{\La-1}_8b_5=0,\allowdisplaybreaks\\&u^{\La-1}_9b_5=\La b_{13},
\quad u^{\La-1}_{10}b_5=(\La-1)b_{15},\quad u^{\La-1}_{11}b_5=0,\allowdisplaybreaks\\
&u^{\La-1}_{12}b_5=\La(\La b_{16}+\pa b_7),\quad
u^{\La-1}_{13}b_5=0,\quad u^{\La-1}_{14}b_5=0,\quad
u^{\La-1}_{15}b_5=0,\allowdisplaybreaks\\
&u^{\La-1}_{16}b_5=\pa b_{15}.
\end{align*}

It follows that in the case $\La\ge 2$ that $N^{E_0,\ov{E}_0}$ is
generated over $\C[\pa]$ by the set
$$S^\La=\{b_5,b_7+b_{10},b_8,b_{11},b_{13},b_{14},b_{15},\La
b_{16}+\pa b_7\}.$$ Hence in this case
$\{b_1,b_2,b_3,b_4,b_6,b_9,b_{10}+2\La\pa b_1,b_{12}\}$ is a
$\C[\pa]$-basis for the $(\C E_0\oplus\C\ov{E}_0)$-invariants of
$M_{\SN^4_+}(-\frac{\La+2}{2},\La,\La)/N$.

The $(L_0,H_0,\ov{H}_0)$-weights of $b_1$, $b_2$, $b_3$, $b_4$,
$b_6$, $b_9$, $b_{10}+2\La\pa b_1$, $b_{12}$ are
$(\Delta,\La,\La)$, $(\Delta+\half,\La+1,\La+1)$,
$(\Delta+\half,\La-1,\La+1)$, $(\Delta+\half,\La+1,\La-1)$,
$(\Delta+1,\La+2,\La)$, $(\Delta+1,\La,\La+2)$,
$(\Delta+1,\La,\La)$, $(\Delta+\frac{3}{2},\La+1,\La+1)$,
respectively.  Hence $M_{\SN^4_+}(-\frac{\La+2}{2},\La)/N$ is a
free $\C[\pa]$-module of rank $8(\La+1)(\La+2)$. So we need to
show that $M_{\SN^4_+}(\frac{\La}{2},\La,\La)/N$ is irreducible.

Again we will study the $(\V_+\oplus sl_2\oplus\ov{sl}_2)$-module
structure of $M_{\NS^4_+}(\frac{\La}{2},\La)/N$.  We can check
directly that $L_n$, for $n\ge 1$, annihilates the vectors $b_1$,
$b_2$, $b_3$, $b_4$, $b_6$, $b_9$, $b_{10}+2\La\pa b_1$,
$b_{12}$. Thus $M_{\SN^4_+}(-\frac{\La+2}{2},\La,\La)/N$ as a
$(\V_+\oplus sl_2\oplus\ov{sl}_2)$-module is a direct sum of the
following eight irreducible modules: $\C[\pa]V_1\cong
L_{\V_+}(-\frac{\La+2}{2})\boxtimes U^{\La,\La}$, $\C[\pa]V_2\cong
L_{\V_+}(-\frac{\La+1}{2})\boxtimes U^{\La+1,\La+1}$,
$\C[\pa]V_3\cong L_{\V_+}(-\frac{\La+1}{2})\boxtimes
U^{\La-1,\La+1}$, $\C[\pa]V_{4}\cong
L_{\V_+}(-\frac{\La+1}{2})\boxtimes U^{\La+1,\La-1}$,
$\C[\pa]V_6\cong L_{\V_+}(-\frac{\La}{2})\boxtimes U^{\La+2,\La}$,
$\C[\pa]V_{9}\cong L_{\V_+}(-\frac{\La}{2})\boxtimes
U^{\La,\La+2}$, $\C[\pa]V_{10}\cong
L_{\V_+}(-\frac{\La}{2})\boxtimes U^{\La,\La}$,
$\C[\pa]V_{12}\cong L_{\V_+}(-\frac{\La-1}{2})\boxtimes
U^{\La+1,\La+1}$, where $V_i$ is the irreducible
$sl_2\oplus\ov{sl}_2$-module generated by $b_i$, for $i\not=10$,
and $V_{10}$ is generated by $b_{10}+2\La\pa b_1 $, and
$U^{\mu,\mu'}$ is the irreducible $sl_2\oplus\ov{sl}_2$-module of
highest weight $(\mu,\mu')$. Note these modules are all
irreducible.  Note further that they are all non-isomorphic.  So
as before to show that $M_{\SN^4_+}(-\frac{\La+2}{2},\La,\La)/N$
is irreducible, it suffices to show that one may send a
$(\V_+\oplus sl_2\oplus\ov{sl}_2)$-highest weight vector in any
irreducible $(\V_+\oplus sl_2\oplus\ov{sl}_2)$-component to the
irreducible component containing the $\G$-highest weight
vectors.  For this purpose we compute
\begin{align*}
&G^{--}_{\half}b_2=2(\La+1)b_1,\quad
\ov{G}^{+-}_{\half}b_3=-2\La(\La+1)b_1,\\
&G^{-+}_{\half}b_4=-2\La(\La+1)b_1,\quad
F_1b_{6}=-2(\La+1)b_1,\\
&\ov{F}_1b_{9}=-2(\La+1)b_1\quad
\ov{F}_1(b_{10}+2\La\pa b_1)=2\La\ov{F}_0 b_1,\\
&\ov{G}^{--}_{\frac{3}{2}}b_{12}=8(\La+1)b_1.
\end{align*}
This settles the case when $\La\ge 2$.

In the case when $\La=1$ $N^{E_0,\ov{E}_0}$ is generated over
$\C[\pa]$ by $$S^\La=\{b_5,b_7+b_{10},b_{13},b_{14},b_{16}+\pa
b_7,\pa b_{15}\}.$$ Therefore
$M_{\SN^4_+}(-\frac{\La+2}{2},\La,\La)/N$ contains a
$\pa$-invariant (and hence $\G$-invariant) vector $b_{15}$. Since
in this case the vectors $b_8=b_{11}=0$,
$M_{\SN^4_+}(-\frac{\La+2}{2},\La,\La)/(N+\C b_{15})$ as a
$\V_+\oplus sl_2\oplus\ov{sl}_2$-module is isomorphic to
$\C[\pa]V_1\oplus\C[\pa]V_2\oplus
\C[\pa]V_3\oplus\C[\pa]V_4\oplus\C[\pa]V_{6}\oplus\C[\pa]V_{9}\oplus
\C[\pa]V_{10}\oplus\C[\pa]V_{12}$.  Every component is
irreducible except for $\C[\pa]V_{12}$, which contains a unique
(irreducible) $\V_+\oplus sl_2\oplus\ov{sl}_2$-submodule
isomorphic to $L_{\V_+}(1)\otimes U^{2,2}$ generated by the
highest weight vector $\pa b_{15}$. But then the above
calculation plus the fact that $$\ov{G}^{--}_{\frac{5}{2}}\pa
b_{12}=24(\La+1)\pa b_1$$ also shows that
$M_{\SN^4_+}(-\frac{\La+2}{2},\La,\La)/(N+\C b_{15})$ is
irreducible.
\end{proof}

We summarize the results of this section in the following theorem.

\begin{thm}
The modules $L_{\SN^4_+}(\Delta,\La,\ov{\La})$, for $\Delta\in\C$
and $\La,\ov{\La}\in\Z_+$, form a complete list of non-isomorphic
finite (over $\C[\pa]$) irreducible $K(1,4)_+$-modules.
Furthermore $L_{\SN^4_+}(\Delta,\La,\ov{\La})$ as a
$\C[\pa]$-module has rank
\begin{itemize}
\item[i.] $8\La(\La+1)$, in the case $2\Delta-\La=0$ and $\La=\ov{\La}$,
\item[ii.] $8(\La+1)(\La+2)$, in the case $2\Delta+\La+2=0$ and $\La=\ov{\La}$.
\item[iii.] $16(\La+1)(\ov{\La}+1)$, in all other cases.
\end{itemize}
Furthermore the $\C[\pa]$-rank of
$L_{\SN^4_+}(\Delta,\La,\ov{\La})_{\bar{0}}$ equals the
$\C[\pa]$-rank of $L_{\SN^4_+}(\Delta,\La,\ov{\La})_{\bar{1}}$ in
all cases.
\end{thm}

\begin{rem}
Again the translation into the languages of modules over conformal
algebras and of conformal modules is straightforward and hence is
omitted. We thus obtain that all finite irreducible modules over
the ``big'' $N=4$ conformal superalgebra are of the form
$L_{\SN^4}(\alpha,\Delta,\La,\ov{\La})$, where
$\alpha,\Delta\in\C$ and $\La,\ov{\La}\in\Z_+$.  Again the
definition of these modules and the action of the conformal
superalgebra on them are easily derived from our explicit
description of a $\C[\pa]$-basis in this section. We note that
the adjoint module is isomorphic to $M_{\SN^4}(0,0,0,0)$.  This
module is not simple, since $K(1,4)$ is not a simple Lie
superalgebra.  Its derived algebra $K(1,4)'$ (which is a simple
formal distribution Lie superalgebra) is an ideal in $K(1,4)$ of
codimension $1$ \cite{KL}. Thus the annihilation subalgebra of
$K(1,4)'$ and $K(1,4)$ are identical, and hence their conformal
modules are identical. Therefore the results in this section also
give explicit description of irreducible conformal modules over
$K(1,4)'$.  We finally remark that the $K(1,4)'$ as a conformal
module over $K(1,4)$ corresponds to $L_{\SN^4}(0,\half,1,1)$.
\end{rem}

\end{document}